\documentclass[12pt]{amsart}

\usepackage{amsmath,amsthm,amssymb,mathtools}
\usepackage{xcolor}
\usepackage{url}
\usepackage{setspace}
\allowdisplaybreaks

\newtheorem{theorem}{Theorem}[section]
\newtheorem{proposition}[theorem]{Proposition}
\newtheorem{lemma}[theorem]{Lemma}
\newtheorem{corollary}[theorem]{Corollary}

\newtheorem{question}[theorem]{Problem}

\theoremstyle{definition}

\newtheorem*{proposition*}{Proposition}
\newtheorem*{theorem*}{Theorem}
\newtheorem*{example*}{Example}

\theoremstyle{remark}
\newtheorem{remark}[theorem]{Remark}
\numberwithin{equation}{section}

\newcommand{\R}{{\mathbb R}}
\newcommand{\Z}{{\mathbb Z}}
\newcommand{\N}{{\mathbb N}}

\newcommand{\x}{{\boldsymbol x}}
\newcommand{\y}{{\boldsymbol y}}

\newcommand{\Wt}{\mathcal{W}}
\newcommand{\Tev}{T_\ast^{\mathrm{ev}}}

\title{Weighted designs of harmonic index}

\author[M.~Hirao]{Masatake Hirao}
\address[M.~Hirao]{Department of Information Science and Technology\\
	Tokyo University of Science\\
	2641 Yamazaki, Noda-shi, Chiba 278-8510\\
	Japan}
\email{hirao@rs.tus.ac.jp}

\author[H.~Nozaki]{Hiroshi Nozaki}
\address[H.~Nozaki]{Department of Mathematics Education\\ 
	Aichi University of Education\\
	1 Hirosawa, Igaya-cho, Kariya, Aichi 448-8542\\
	Japan}
\email{hnozaki@auecc.aichi-edu.ac.jp}

\author[K.~Tasaka]{Koji Tasaka}
\address[K.~Tasaka]{Department of Mathematics\\
	Kindai University\\
	Osaka 577-8502\\
	Japan\\
	Research Institute for Mathematical Sciences\\
	Kyoto University\\
	Kyoto 606-8502\\ 
	Japan}
\email{tasaka@math.kindai.ac.jp,\ tasaka@kurims.kyoto-u.ac.jp}

\author[S.~Wakatsuki]{Satoshi Wakatsuki}
\address[S.~Wakatsuki]{Faculty of Mathematics and Physics\\
Institute of Science and Engineering\\
Kanazawa University\\
Kakumamachi, Kanazawa, Ishikawa, 920-1192\\ 
Japan}
\email{wakatsuk@staff.kanazawa-u.ac.jp}


\begin{document}

\begin{abstract}
Let $X$ be a finite subset of the unit sphere $\mathbb{S}^{d-1}$, $d\ge3$.
We introduce the \emph{weighted harmonic strength} $T_\ast(X)$: the set of positive integers $\ell$ such that, for some positive weight $w$ on $X$, the weighted sum $\sum_{\x\in X}w(\x)P(\x)$ vanishes for every harmonic polynomial $P$ of degree $\ell$; here the weight may depend on $\ell$.
It contains the harmonic strength $T(X)$, and the two sets agree whenever $X$ is a single orbit of a finite subgroup of the orthogonal group.
Our main tool is a criterion, uniform in the weight, which excludes from $T_\ast(X)$ every even index beyond an explicit threshold; the threshold is obtained from the Laplace integral representation of the Gegenbauer polynomials, and can be determined in exact arithmetic.
Together with a reduction of the design condition to a Gram matrix attached to the orbits of a symmetry group, this turns the determination of $T_\ast(X)$ into a finite computation.
We carry it out for unions of normalized shells of $\Z^3$, of $D_4\cup\sqrt{2}D_4^\ast$ and of $E_8$.
In particular, the point sets underlying the spherical cubature formulas of de la Harpe, Pache and Venkov carry a long list of sporadic harmonic indices far beyond their design strength.
\end{abstract}

\maketitle

\section{Introduction}\label{sec:intro}

\subsection{The problem}\label{subsec:intro-problem}
Throughout the paper we assume $d\ge3$.
We denote by $\N=\{1,2,3,\ldots\}$ the set of positive integers, by $2\N$ the subset of even ones and by $\mathbb{O}:=\N\setminus2\N$ the subset of odd ones.
A cubature formula on the unit sphere $\mathbb{S}^{d-1}:=\{\x=(x_1,\ldots,x_d)\in\R^d\mid x_1^2+\cdots+x_d^2=1\}$ of the $d$-dimensional Euclidean space $\R^d$ replaces the normalized surface integral of a polynomial by a finite weighted sum of its values.
Let $\sigma$ be the surface measure on $\mathbb{S}^{d-1}$ and $|\mathbb{S}^{d-1}|=\int_{\x \in \mathbb{S}^{d-1}}d\sigma(\x)$.
Following the terminology of algebraic combinatorics, we call a finite subset $X\subset\mathbb{S}^{d-1}$ a \emph{weighted $t$-design} if there is a positive weight $w:X\rightarrow \R_{>0}$ on $X$ with $\sum_{\x\in X}w(\x)=1$ for which
\begin{equation}\label{eq:weighted_design}
 \sum_{\x\in X} w(\x) f(\x) = \frac{1}{|\mathbb{S}^{d-1}|} \int_{\x \in \mathbb{S}^{d-1}}f(\x)\,d\sigma(\x) 
\end{equation}
holds for all polynomials $f\in \R[x_1,\ldots,x_d]$ of degree $\le t$.
If \eqref{eq:weighted_design} with the equal weight $w\equiv1/|X|$ holds for all polynomials $f\in \R[x_1,\ldots,x_d]$ of degree $\le t$, then $X$ is called a \emph{spherical $t$-design}. 
Spherical designs were introduced by Delsarte--Goethals--Seidel~\cite{DGS77}; cubature formulas on the sphere go back further in numerical analysis, at least to Sobolev~\cite{Sob62}.
One of the basic questions about weighted designs is the determination of the minimum cardinality required for a weighted $t$-design of given $t$.
Lower bounds for this cardinality were obtained by M\"{o}ller~\cite{Mol76,Mol79}; the equal-weight case is due to the independent work of Delsarte--Goethals--Seidel~\cite{DGS77}.
For the history and theory of weighted designs (or cubature formulas on the sphere), see \cite{BB09,BBHS10} and the references therein.
 
We begin with an equivalent description of weighted $t$-designs in terms of harmonic polynomials.
Condition~\eqref{eq:weighted_design} is equivalent to the vanishing of the weighted sums
\begin{equation}\label{eq:intro-harm}
 \sum_{\x\in X}w(\x)P(\x)=0,
 \qquad \forall P\in\mathrm{Harm}_\ell,
\end{equation}
for every $\ell$ with $1\le\ell\le t$, where $\mathrm{Harm}_\ell$ denotes the space of real homogeneous harmonic polynomials of degree $\ell$ in $d$ variables
\[ \mathrm{Harm}_\ell:=\left\{P\in \R[x_1,\ldots,x_d]\ \middle|\ \sum_{j=1}^d \frac{\partial^2}{\partial x_j^2}P(\x)=0,\ P(\lambda \x)=\lambda^\ell P(\x)\right\},\]
regarded as functions on the sphere $\mathbb{S}^{d-1}$.
See Proposition \ref{prop:weighted_design_equiv}.
 
In this paper we are interested in determining all degrees $\ell>t$ at which \eqref{eq:intro-harm} continues to hold.
More precisely, for a finite subset $X\subset\mathbb{S}^{d-1}$, let
\[
 \Wt(X):= \left\{ w:X\rightarrow \R_{>0}\ \middle| \ \sum_{\x\in X} w(\x)=1\right\}
\] 
be the set of positive weights on $X$, and consider the two subsets of $\N$
\begin{equation}\label{eq:strength}
\begin{aligned}
 T(X)&:=\bigl\{\ell\in\N\bigm|\ \text{\eqref{eq:intro-harm} holds for }
 w\equiv1/|X|\bigr\},\\
 T_\ast(X)&:=\bigl\{\ell\in\N\bigm|\ \text{\eqref{eq:intro-harm} holds for
 some }w\in\Wt(X)\bigr\}.
\end{aligned}
\end{equation}
We call $T(X)$ the \emph{harmonic strength} of $X$ and $T_\ast(X)$ its \emph{weighted harmonic strength}.
We call $\ell\in\N$ an \emph{index} of $X$ if $\ell\in T_\ast(X)$, and we say that $w\in\Wt(X)$ \emph{realizes} $\ell$ if \eqref{eq:intro-harm} holds for this $w$.
In $T_\ast(X)$ the weight is allowed to depend on the index $\ell$; to the best of our knowledge, this setting has not been considered before.
The problem we address is the following.
 
\begin{question}\label{q:main}
Given a finite subset $X\subset \mathbb{S}^{d-1}$, determine $T(X)$ and $T_\ast(X)$.
\end{question}
 
Clearly $T(X)\subset T_\ast(X)$. Moreover, $T(gX)=T(X)$ and $T_\ast(gX)=T_\ast(X)$ for every $g$ in the orthogonal group $O(\R^d)$.
As a starting point, we show in Proposition~\ref{prop:averaging} that $T(X)=T_\ast(X)$ holds whenever a finite subgroup of $O(\R^d)$ acts transitively on $X$.
Hence a gap between $T(X)$ and $T_\ast(X)$, the situation we are interested in, can occur only for finite subsets $X$ that are not single orbits.
Numerous examples of such subsets are provided by unions of normalized shells of a Euclidean lattice.
Goethals--Seidel~\cite{GS81} and de la Harpe, Pache and Venkov~\cite{HPV07} initiated the construction of cubature formulas supported on such unions.
In this paper, we revisit their constructions as test cases for Problem~\ref{q:main}.
 
As background, the notion of harmonic strength already appears in the study of \emph{spherical $T$-designs}, which for a given $T\subset\N$ are the finite subsets $X$ with $T\subset T(X)$.
These were introduced by Delsarte--Seidel~\cite{DS89} as a spherical analogue of
the $T$-designs in association schemes \cite[Section~3.4]{Del73}; their systematic study was begun by Bannai--Okuda--Tagami~\cite{BOT15}; see also \cite{HNT23,Mie13,OY16,Pan22,ZBBKY17}.
The version with positive weights, called a \emph{weighted $T$-design}, was considered by Andreev~\cite{And00} and, for association schemes, by Delsarte--Levenshtein~\cite{DL98} under the name of $\tau$-designs.
Note that the condition $T\subset T_\ast(X)$ is weaker than the condition that $X$ be a weighted $T$-design, since the latter requires a single $w$ to work for all $\ell\in T$.
 
Problem~\ref{q:main} for $T(X)$ originates from the work of Boyvalenkov, Danev and Kazakov~\cite{BDK01}, who determined $T(X)$ for several classical examples, including the vertex sets of the icosahedron and of the regular $600$-cell and the minimal vectors of the Leech lattice.
For example, for the vertex set $X$ of the $600$-cell, one has $T(X)\cap2\N=\{2,4,6,8,10,14,16,18,22,26,28,34,38,46,58\}$.
The complete determination of $T(X)$ for finite sets $X \subset \mathbb{S}^3$ arising from exceptional quaternionic unit groups, including the $D_4$ root system, is carried out in~\cite{HNT}.

\subsection{Main results}\label{subsec:intro-results}
If a finite subset $X\subset\mathbb{S}^{d-1}$ is \emph{antipodal}, that is, $-X=X$ where $-X:=\{-\x\mid\x\in X\}$, then it follows that $\mathbb{O}\subset T(X)$.
Hence,
\[ \mathbb{O}\subset T_\ast(X)\quad \mbox{if $X$ is antipodal}.\]
Since we are considering unions of shells of a lattice, which are antipodal, only the even elements of $T(X)$ and $T_\ast(X)$ remain to be determined, and we let
\[ \Tev(X):=T_\ast(X)\cap 2\N .\]

Our main result provides an explicit $\ell_0$ such that $2\ell\not\in\Tev(X)$ for every $\ell\ge\ell_0$.
In particular, we have $|\Tev(X)|<\infty$.
We remark that, in the equal-weight setting, the finiteness of $T(X)\cap2\N$ was essentially shown by Boyvalenkov, Danev and Kazakov~\cite{BDK01}.

Let us denote the inner product set of $X$ by
\begin{equation}\label{eq:A(X)}
A(X):=\{\langle\x,\y\rangle\mid\x,\y\in X,\ \x\ne\y\},
\end{equation}
where $\langle \x,\y\rangle$ denotes the Euclidean inner product of vectors $\x,\y\in \R^d$, and, for $\ell\in\N$ and $a\in[0,1)$, put
\begin{equation}\label{eq:intro-FLap}
 F^{\mathrm{Lap}}_{2\ell,d}(a)
 :=\frac1{s_d}\int_0^\pi\bigl(a^2+(1-a^2)\cos^2\theta\bigr)^\ell
 \sin^{d-3}\theta\,d\theta,
 \qquad
 s_d:=\int_0^\pi\sin^{d-3}\theta\,d\theta .
\end{equation}
The function $F^{\mathrm{Lap}}_{2\ell,d}(a)$ comes from the Laplace integral representation of the Gegenbauer polynomial.
A crucial feature is that, for each $a\in[0,1)$, this function satisfies $F^{\mathrm{Lap}}_{2(\ell+1),d}(a)\le F^{\mathrm{Lap}}_{2\ell,d}(a)$ for all $\ell\ge1$, and $F^{\mathrm{Lap}}_{2\ell,d}(a)\to0$ as $\ell\to\infty$.
This allows us to find $\ell_0$ as in the following theorem.

\begin{theorem}\label{thm:main}
Let $d\ge3$ and $X\subset\mathbb{S}^{d-1}$ be a finite antipodal subset.
Then $\Tev(X)$ is finite.
Moreover, suppose that $|X|\ge4$, and let $a_X$ be the largest $|a|$ with $a\in A(X)\setminus\{-1\}$.
If $\ell_0\in \N$ denotes the least positive integer satisfying
\[
 F^{\mathrm{Lap}}_{2\ell_0,d}\bigl(a_X\bigr)<\frac{2}{|X|-2},
\]
then we have
\[
 \Tev(X)\subset\{2\ell\mid \ell\in\N,\ \ell<\ell_0\}.
\]
\end{theorem}

We call $\ell_0$ in Theorem \ref{thm:main} the \emph{Laplace cutoff} of $X$.
Its size is governed by the asymptotics of \eqref{eq:intro-FLap}: applying Laplace's method at the two endpoints $\theta=0,\pi$ gives, for each fixed $a\in[0,1)$,
\begin{equation*}
 F^{\mathrm{Lap}}_{2\ell,d}(a)
 \sim\frac{\Gamma\!\left(\frac{d-1}{2}\right)}{\sqrt{\pi}}
 \bigl(\ell\,(1-a^2)\bigr)^{-\frac{d-2}{2}}
 \qquad(\ell\to\infty),
\end{equation*}
so that the Laplace cutoff of $X$ is of the expected size
\begin{equation*}
 \ell_0\ \approx\ \frac{1}{1-a_X^2}
 \left(\frac{\Gamma\!\left(\frac{d-1}{2}\right)\bigl(|X|-2\bigr)}{2\sqrt{\pi}}\right)^{\frac{2}{d-2}} .
\end{equation*}

Theorem~\ref{thm:main} leaves the finitely many candidates $2,4,\ldots,2\ell_0-2$, and it remains to decide which of them are indices of $X$. 
We carry this out for unions of lattice shells.
For a Euclidean lattice $\Lambda\subset\R^d$ and $n>0$ write
\begin{equation}\label{eq:shell}
 \Lambda[n]:=\frac1{\sqrt n}\{\x\in\Lambda\mid\langle\x,\x\rangle=n\},
 \qquad
 \Lambda[n_1,\ldots,n_r]:=\bigcup_{j=1}^r\Lambda[n_j],
\end{equation}
and use the same notation when $\Lambda$ is a finite union of lattices.
 
The smallest example already exhibits the two phenomena studied in this paper: a union of shells can have more indices than each of its shells, and the new indices need not all be realized by the same weight.
The octahedron $\Z^3[1]$, the cube $\Z^3[3]$ and their union all have harmonic strength $\{2\}\cup\mathbb{O}$, whereas
\[
 T_\ast\bigl(\Z^3[1,3]\bigr)=\{2,4,10\}\cup\mathbb{O},
\]
and no single weight realizes both of the new indices $4$ and $10$ (Theorem~\ref{thm:Z3-shells}).
 
Our main examples are the point sets underlying the cubature formulas of~\cite{HPV07}, which we treat in Appendix~\ref{app:HPV}.
For $\Lambda=D_4\cup\sqrt{2}D_4^\ast$, the sets $\Lambda[2]$, $\Lambda[2,6]$ and $\Lambda[2,6,10]$ are, by \cite[Theorem~$D_4$]{HPV07}, a spherical $7$-design and weighted designs of strengths $11$ and $15$, respectively.
All three have the same harmonic strength $\{2,4,6,10,14,22\}\cup\mathbb{O}$, whereas their weighted harmonic strengths grow with the number of shells; the largest gains the twelve further even indices $8,12,18,26,30,34,38,42,46,54,62,70$ (Theorem~\ref{thm:D4-shells}).
For $\Lambda=E_8$, the three sets $\Lambda[2,4]$, $\Lambda[2,6,8]$ and $\Lambda[2,4,6,8]$ of \cite[Theorem~$E_8$]{HPV07} all have harmonic strength $\{2,4,6,10\}\cup\mathbb{O}$, and the last of them has twelve further even indices (Theorem~\ref{thm:E8-shells}).
The Laplace cutoffs of these six sets range from $24$ to $6600$, so the statements rest on deciding several thousand candidate indices in exact arithmetic.

\section{Comparison of the two strengths}\label{sec:orbits}

In this section, we show that the two strengths $T(X)$ and $T_\ast(X)$ defined in \eqref{eq:strength} coincide when a finite subgroup of the orthogonal group acts transitively on $X\subset \mathbb{S}^{d-1}$.

For this, let $O(\R^d)$ denote the group of orthogonal transformations of $\R^d$.
The group $O(\R^d)$ acts on the sphere $\mathbb{S}^{d-1}$ by $(g,\x)\mapsto g(\x)$.
This induces an action of the group $O(\R^d)$ on the space $\mathrm{Harm}_\ell$ of harmonic polynomials.
More precisely, for $P\in\mathrm{Harm}_\ell$ and $g\in O(\R^d)$, define the polynomial $gP$ by $(gP)(\x)=P(g^{-1}(\x))$.
Then, $gP$ again lies in $\mathrm{Harm}_\ell$.

\begin{proposition}\label{prop:averaging}
Let $X\subset \mathbb{S}^{d-1}$ be a finite subset, let $H\subset O(\R^d)$ be a finite subgroup preserving $X$, and let $X=\mathcal{O}_1\sqcup\cdots\sqcup\mathcal{O}_r$ be the decomposition of $X$ into $H$-orbits.
Then, for $\ell\in\N$, we have $\ell\in T_\ast(X)$ if and only if there are positive real numbers $c_1,\ldots,c_r$ with $\sum_{i=1}^r c_i=1$ such that
\begin{equation}\label{eq:orbit_condition}
 \sum_{i=1}^r\frac{c_i}{|\mathcal{O}_i|}\sum_{\x\in\mathcal{O}_i}P(\x)=0
\end{equation}
for every $P\in\mathrm{Harm}_\ell$.
In particular, if $H$ acts transitively on $X$, then
\[
 T_\ast(X)=T(X).
\]
\end{proposition}

\begin{proof}
Suppose first that such $c_1,\ldots,c_r$ exist, and define $w:X\to\R_{>0}$ by $w(\x):=c_i/|\mathcal{O}_i|$ for $\x\in\mathcal{O}_i$.
Then $w(\x)>0$ and $\sum_{\x\in X}w(\x)=\sum_{i=1}^r c_i=1$, so $w\in\Wt(X)$, and $\sum_{\x\in X}w(\x)P(\x)$ is exactly the left-hand side of \eqref{eq:orbit_condition}.
Hence $\ell\in T_\ast(X)$.
 
Conversely, let $\ell\in T_\ast(X)$ and $w\in\Wt(X)$ realize $\ell$.
Define the $H$-average
\[
 w^{H}(\x):=\frac{1}{|H|}\sum_{h\in H}w(h(\x))\qquad(\x\in X),
\]
which defines a map $w^H:X\rightarrow \R_{>0}$.
Since $\x\mapsto h(\x)$ is a bijection of $X$ for each $h\in H$, we obtain
\[
 \sum_{\x\in X}w^{H}(\x)
 =\frac{1}{|H|}\sum_{h\in H}\sum_{\x\in X}w(h(\x))
 =\frac{1}{|H|}\sum_{h\in H}\sum_{\y\in X}w(\y)=1 ,
\]
so $w^{H}\in\Wt(X)$.
By construction $w^{H}(h(\x))=w^{H}(\x)$ for all $h\in H$, that is, $w^{H}$ is constant on each orbit $\mathcal{O}_i$.

Using the weight $w^{H}$, we define $c_i=\sum_{\x\in\mathcal{O}_i}w^{H}(\x)$ for $i=1,2,\ldots,r$.
Then, $c_i>0$ and $\sum_{i=1}^r c_i=1$.
It follows that $w^H(\x)=c_i/|\mathcal{O}_i|$ if $\x\in\mathcal{O}_i$.
Hence, for $ P\in \mathrm{Harm}_\ell$, one has
\begin{align*}
\sum_{i=1}^r\frac{c_i}{|\mathcal{O}_i|}\sum_{\x\in\mathcal{O}_i}P(\x) &= \sum_{\x\in X} w^H(\x)P(\x)\\
& =\frac{1}{|H|}\sum_{\x\in X}\sum_{h\in H}w(h(\x))P(\x)\\
 &=\frac{1}{|H|}\sum_{h\in H}\sum_{\y\in X}w(\y)\,(hP)(\y).
\end{align*}
Since $\ell \in T_\ast(X)$ with the weight $w$ and $hP\in \mathrm{Harm}_\ell$, the sum $\sum_{\y\in X}w(\y)\,(hP)(\y)$ in the last term is 0 for each $h\in H$,
which proves \eqref{eq:orbit_condition}.
 
Finally, if $H$ acts transitively on $X$, then $r=1$ and $c_1=1$, so \eqref{eq:orbit_condition} reads
\[
 \frac{1}{|X|}\sum_{\x\in X}P(\x)=0
 \qquad\text{for every }P\in\mathrm{Harm}_\ell,
\]
which is precisely the condition $\ell\in T(X)$.
Hence $T_\ast(X)\subset T(X)$, and the reverse inclusion always holds.
\end{proof}

Proposition \ref{prop:averaging} provides a large supply of finite sets with $T(X)=T_\ast(X)$: every orbit of a finite subgroup of $O(\R^d)$ is such a set.
This covers most of the finite subsets for which $T(X)$ has been determined in the literature, among them the vertex sets of the icosahedron and of the $600$-cell, the shortest vectors of the Leech lattice \cite{BDK01}, and those arising from the finite quaternionic unit groups treated in \cite{HNT}.
In all of these cases, relaxing the condition by allowing arbitrary positive weights does not enlarge the harmonic strength.
Taking $X$ to be a union of normalized shells of a lattice, we examine whether the inclusion $T(X)\subset T_\ast(X)$ is strict.

\begin{remark}\label{rem:euclidean}
We note that Proposition~\ref{prop:averaging} has an obvious analogue for Euclidean designs introduced in \cite{NS88}.
Weighted designs or Euclidean designs on which a finite subgroup $G$ of $O(\R^d)$ acts have often been constructed by taking $X$ to be a union of $G$-orbits and assigning the same weight to all points of each orbit (see, e.g., \cite{Sob62,GS81,HPV07,NS12}); in \cite[Section~4]{NS12}, for instance, this is assumed throughout.
Proposition~\ref{prop:averaging} shows that this assumption is no loss of generality: if $X$ is $G$-invariant and a positive weight realizes $\ell$, then so does a positive weight that is constant on each $G$-orbit.
Although the proof is a simple averaging argument, this fact does not seem to
be widely known, and we have not found it stated explicitly in the literature.
\end{remark}

\section{Weighted harmonic strength via Gegenbauer polynomials}\label{sec:gegen}
In this section, we give another description of the strengths $T(X)$ and $T_\ast(X)$ in terms of Gegenbauer polynomials. 
This description is used to prove Theorem \ref{thm:main}.

For $\ell\in \Z_{\ge0}$, the \emph{Gegenbauer polynomial} $G_{\ell,d}(x)\in \R[x]$ is defined by the three-term recurrence
 \begin{equation*}
 (\ell+d-2)\,G_{\ell+1,d}(x)=(2\ell+d-2)\,x\,G_{\ell,d}(x)-\ell\,G_{\ell-1,d}(x)
 \qquad(\ell\ge1)
\end{equation*}
with the initial polynomials $G_{0,d}(x)=1$ and $G_{1,d}(x)= x$.
For example, 
\begin{align*}
&G_{2,d}(x)=\frac{d x^2-1}{d-1},\quad G_{3,d}(x)=\frac{x \left((d+2) x^2-3\right)}{d-1},\\
& G_{4,d}(x)=\frac{(d+2) x^2 \left((d+4) x^2-6\right)+3}{d^2-1}.
\end{align*}
In terms of the Gegenbauer polynomial $C_\ell^{(\lambda)}(x)$ in the classical normalization, we have
\[
 G_{\ell,d}(x)=\frac{C_\ell^{((d-2)/2)}(x)}{C_\ell^{((d-2)/2)}(1)}.
\]
Alternatively, for an integer $\ell\ge0$, let $P_\ell^{(\alpha,\beta)}(x)$ be the Jacobi polynomial of degree $\ell$ (see Szeg\H{o}'s book \cite{Sze75} for the definition) with the normalization $P_\ell^{(\alpha,\alpha)}(1)=\binom{\ell+\alpha}{\ell}$.
Then we have
\[
 G_{\ell,d}(x) = \frac{P_\ell^{((d-3)/2,(d-3)/2)}(x)}{\binom{\ell+(d-3)/2}{\ell}}.
\]
With our normalization, we have $G_{\ell,d}(1)=1$ and $G_{\ell,d}(-1)=(-1)^\ell$ for all $\ell\ge0$.
The polynomial $G_{\ell,d}(x)$ is an even (resp.~odd) polynomial of degree $\ell$ if $\ell$ is even (resp.~odd).

The following characterizations of weighted $t$-designs are well known.

\begin{proposition}
\label{prop:weighted_design_equiv}
Let $X\subset \mathbb{S}^{d-1}$ be a finite subset, let $t\in\N$, and let $w\in \Wt(X)$.
Then the following conditions are equivalent{\rm :}
\begin{itemize}
\item[$(\mathrm{i})$] the equality \eqref{eq:weighted_design} holds for all $f\in \R[x_1,\ldots,x_d]$ of degree $\le t${\rm ;} \item[$(\mathrm{ii})$] $\displaystyle \sum_{\x,\y \in X} w(\x)\, w(\y)\,G_{\ell,d}(\langle \x,\y\rangle) = 0$ for all $\ell=1,2,\ldots,t${\rm ;} \item[$(\mathrm{iii})$] $\displaystyle \sum_{\x\in X} w(\x)\,P(\x) = 0$ for all $P\in\mathrm{Harm}_\ell$ and all $\ell=1,2,\ldots,t$.
\end{itemize}
In particular, $X$ is a weighted $t$-design if and only if there exists a weight $w\in\Wt(X)$ satisfying $(\mathrm{ii})$, or equivalently $(\mathrm{iii})$.
\end{proposition}

We briefly sketch the proof of Proposition \ref{prop:weighted_design_equiv}.
The equivalence of $(\mathrm{i})$ and $(\mathrm{iii})$ in Proposition \ref{prop:weighted_design_equiv} goes back to Delsarte--Goethals--Seidel~\cite[Theorem~5.2]{DGS77} in the case of equal weights; for weighted designs, see de la Harpe, Pache and Venkov~\cite[Proposition~4]{HPV07}, and also \cite{NS88}.
Indeed, restricted to $\mathbb{S}^{d-1}$, every $f\in\R[x_1,\ldots,x_d]$ of degree $\le t$ decomposes as $f=\sum_{\ell=0}^{t}P_\ell$ with $P_\ell\in\mathrm{Harm}_\ell$; since $\sum_{\x\in X}w(\x)=1$ and $\int_{\mathbb{S}^{d-1}}P_\ell\,d\sigma=0$ for $\ell\ge1$, one has
\[
 \sum_{\x\in X} w(\x) f(\x)
 - \frac{1}{|\mathbb{S}^{d-1}|} \int_{\x \in \mathbb{S}^{d-1}}f(\x)\,d\sigma(\x)
 =\sum_{\ell=1}^{t}\ \sum_{\x\in X} w(\x) P_\ell(\x),
\]
and $(\mathrm{iii})\Rightarrow(\mathrm{i})$ follows.
The implication $(\mathrm{i})\Rightarrow(\mathrm{iii})$ is immediate from $\int_{\mathbb{S}^{d-1}}P_\ell\,d\sigma=0$ for $\ell\ge1$.
The equivalence of $(\mathrm{ii})$ and $(\mathrm{iii})$ is a consequence of the addition formula for spherical harmonics; we record it as Lemma~\ref{lem:harmonic} below, since it will be used repeatedly.
For equal weights, condition $(\mathrm{ii})$ appears in \cite[Theorem~5.5]{DGS77}.

We equip $\mathrm{Harm}_\ell$ with the inner product
\begin{equation}\label{eq:harm-inner}
 (P,Q):=\frac{1}{|\mathbb{S}^{d-1}|}\int_{\x\in\mathbb{S}^{d-1}}P(\x)Q(\x)\,d\sigma(\x),
\end{equation}
and we fix an orthonormal basis $\varphi_{\ell,1},\ldots,\varphi_{\ell,h_\ell}$ of $\mathrm{Harm}_\ell$, where $h_\ell:=\dim\mathrm{Harm}_\ell$.
Then the addition formula (see, e.g., \cite[Theorem~3.3]{DGS77}) states that
\begin{equation}\label{eq:addition}
 \sum_{j=1}^{h_\ell}\varphi_{\ell,j}(\x)\,\varphi_{\ell,j}(\y)
 =h_\ell\,G_{\ell,d}(\langle\x,\y\rangle)
 \qquad(\x,\y\in\mathbb{S}^{d-1}).
\end{equation}

\begin{lemma}\label{lem:harmonic}
Let $X\subset\mathbb{S}^{d-1}$ be a finite subset, let $\ell\in\N$, and let $w\in\Wt(X)$.
Then
\begin{equation}\label{eq:sos}
 \sum_{\x,\y\in X}w(\x)\,w(\y)\,G_{\ell,d}(\langle\x,\y\rangle)
 =\frac{1}{h_\ell}\sum_{j=1}^{h_\ell}
 \Bigl(\sum_{\x\in X}w(\x)\,\varphi_{\ell,j}(\x)\Bigr)^{2}\;\ge\;0 ,
\end{equation}
and the left-hand side vanishes if and only if
\[
 \sum_{\x\in X}w(\x)\,P(\x)=0
 \qquad\text{for every }P\in\mathrm{Harm}_\ell .
\]
In particular, the weighted harmonic strength defined in \eqref{eq:strength} is given by
\begin{equation}\label{eq:whs}
  T_\ast (X) = \left\{ \ell \in \N \; \middle| \;
  \exists\, w \in \Wt(X)\; \text{s.t.}\;
  \sum_{\x,\y \in X} w(\x)\,w(\y)\,G_{\ell,d}(\langle \x,\y\rangle) = 0
  \right\}.
\end{equation}
\end{lemma}
 
\begin{proof}
Substituting \eqref{eq:addition} gives
\begin{align*}
 \sum_{\x,\y\in X}w(\x)w(\y)G_{\ell,d}(\langle\x,\y\rangle)
& =\frac{1}{h_\ell}\sum_{j=1}^{h_\ell}\sum_{\x,\y\in X}w(\x)w(\y)
 \varphi_{\ell,j}(\x)\varphi_{\ell,j}(\y)\\
& =\frac{1}{h_\ell}\sum_{j=1}^{h_\ell}\Bigl(\sum_{\x\in X}w(\x)\varphi_{\ell,j}(\x)\Bigr)^{2},
\end{align*}
from which \eqref{eq:sos} follows.
Consequently, the left-hand side of \eqref{eq:sos} vanishes if and only if
\[
 \sum_{\x\in X}w(\x)\varphi_{\ell,j}(\x)=0
 \qquad(j=1,\ldots,h_\ell).
\]
Since $\varphi_{\ell,1},\ldots,\varphi_{\ell,h_\ell}$ is a basis of $\mathrm{Harm}_\ell$, we obtain the desired equivalence.
\end{proof}

\section{Finiteness of even weighted harmonic strengths}\label{sec:upper}

The aim of this section is to present our method for computing, for an antipodal set $X\subset\mathbb{S}^{d-1}$, an integer $\ell_0$ such that $2\ell\notin\Tev(X)$ for all $\ell\ge\ell_0$.
The key ingredients are \eqref{eq:whs} and the fact that $G_{2\ell,d}(a)\to0$ as $\ell\to\infty$ for every $a\in A(X)\setminus\{-1\}$.
To make this computation effective, one needs to choose an explicit family $\{F_{2\ell,d}\}_{\ell\in\N}$ satisfying certain conditions.
In the next section, we construct such a family and prove Theorem~\ref{thm:main}.

In what follows, $X$ denotes a finite antipodal subset of $\mathbb{S}^{d-1}$; this is the only case needed later, and the argument extends to the general case with minor modifications (see Remark~\ref{rem:non_antipodal}).
To cover that case as well, we extend every $w\in\Wt(X)$ to $\mathbb{S}^{d-1}$ by setting $w(\y):=0$ for $\y\notin X$; then the quantity $\Psi_{X,w}$ below, as well as Lemma~\ref{lem:mass}, makes sense for an arbitrary finite subset $X\subset\mathbb{S}^{d-1}$.

Using the set $A(X)$ of inner products of distinct elements of $X$ defined in \eqref{eq:A(X)}, we first write the left-hand side of \eqref{eq:sos} as
\[
  \sum_{\x,\y\in X} w(\x)\,w(\y)\,G_{2\ell,d}(\langle\x,\y\rangle)
  = \sum_{a\in A(X)}\sum_{\substack{\x,\y\in X\\ \langle\x,\y\rangle=a}} w(\x)\,w(\y)\,G_{2\ell,d}(a)
  + \sum_{\x\in X} w(\x)^2\,G_{2\ell,d}(1).
\]
Note that we have $G_{2\ell,d}(1)=G_{2\ell,d}(-1)=1$ for all $\ell\ge0$.
For $a\in A(X)$ and $w\in\Wt(X)$, we define
\[
  \Phi_{X,w}(a) := \sum_{\substack{\x,\y\in X\\ \langle\x,\y\rangle=a}} w(\x)\,w(\y)
\]
and
\[
  \Psi_{X,w} := \sum_{\x\in X} w(\x)^2 + \sum_{\x\in X} w(\x)\,w(-\x).
\]
Since $\langle\x,\y\rangle=-1$ holds for $\x,\y\in\mathbb{S}^{d-1}$ if and only if $\x=-\y$, the second term of $\Psi_{X,w}$ coincides with $\Phi_{X,w}(-1)$.
Thus we obtain
\begin{equation}\label{eq:gegen_psi}
  \sum_{\x,\y\in X} w(\x)\,w(\y)\,G_{2\ell,d}(\langle\x,\y\rangle)
  = \Psi_{X,w} + \sum_{a\in A(X)\setminus\{-1\}} \Phi_{X,w}(a)\,G_{2\ell,d}(a).
\end{equation}
Note that, since $w(\x)>0$ for all $\x\in X$, we have $\Phi_{X,w}(a)>0$ for all $a\in A(X)$ and $\Psi_{X,w}>0$.

We now prove two lemmas.
The first one is an identity relating $\Phi_{X,w}(a)$ and $\Psi_{X,w}$.

\begin{lemma}\label{lem:mass}
Let $X\subset\mathbb{S}^{d-1}$ be a finite subset.
Then, for any $w\in \Wt(X)$, we have
\begin{equation}\label{eq:mass}
  \sum_{a\in A(X)\setminus\{-1\}}\Phi_{X,w}(a) = 1-\Psi_{X,w}.
\end{equation}
In particular $\Psi_{X,w}\le 1$, and equality holds if and only if $A(X)\setminus\{-1\}=\emptyset$.
\end{lemma}
 
\begin{proof}
One computes
\[
  1=\Bigl(\sum_{\x\in X}w(\x)\Bigr)^{2}
  =\sum_{\x,\y\in X}w(\x)w(\y)
  =\Psi_{X,w}+\sum_{a\in A(X)\setminus\{-1\}}\Phi_{X,w}(a),
\]
from which the identity \eqref{eq:mass} follows.
All the terms above are nonnegative, so $\Psi_{X,w}\le1$; moreover, since $\Phi_{X,w}(a)>0$ for every $a\in A(X)$, the sum on the left-hand side of \eqref{eq:mass} vanishes if and only if $A(X)\setminus\{-1\}=\emptyset$.
\end{proof}

The second lemma computes the infimum
\[
 \Psi_X^{\inf}:=\inf_{w\in \Wt(X)}\Psi_{X,w}
\]
for an antipodal set $X$.

\begin{lemma}\label{lem:Psi_inf}
Let $X\subset\mathbb{S}^{d-1}$ be a finite antipodal subset.
Then
\begin{equation}\label{eq:inf_Psi}
 \Psi_X^{\inf}=\frac{2}{|X|},
\end{equation}
and the infimum is attained by the equal weight $w\equiv1/|X|$.
\end{lemma}

\begin{proof}
Let $w\in \Wt(X)$.
Write $X = \{\x_1, -\x_1, \ldots, \x_m, -\x_m\}$, and set $u_i := w(\x_i)+ w(-\x_i)$ for $i = 1, \ldots, m$.
Then $u_i > 0$, $\sum_{i=1}^m u_i = 1$, and
\begin{equation}\label{eq:Psi}
  \Psi_{X,w} = \sum_{i=1}^m u_i^2.
\end{equation}
By the Cauchy--Schwarz inequality,
\[
  \sum_{i=1}^m u_i^2
  \geq \frac{\left(\sum_{i=1}^m u_i\right)^2}{m}
  = \frac{1}{m},
\]
with equality if and only if $u_1 = \cdots = u_m = 1/m$.
This condition is satisfied by the equal weight $w \equiv 1/|X| = 1/(2m)$, for which $u_i = 1/m$ holds for all $i$.
Therefore $\Psi_X^{\inf}=1/m=2/|X|$, attained by $w\equiv1/|X|$.
\end{proof}

We now give a framework for computing an integer $\ell_0$ such that $2\ell\notin\Tev(X)$ for every $\ell\ge\ell_0$.
Let $\mathcal{B}_d$ denote the set of families of functions $\{F_{2\ell,d}\}_{\ell\in\N}$ on $[0,1)$ satisfying the following conditions for all $a,b \in [0,1)$ and $\ell\in\N$:
\begin{itemize}
\item[(B1)] $\displaystyle\lim_{\ell \to \infty} F_{2\ell,d}(a) = 0$, \item[(B2)] $|G_{2\ell,d}(a)| \leq F_{2\ell,d}(a)$, \item[(B3)] $F_{2(\ell+1),d}(a) \leq F_{2\ell,d}(a)$, \item[(B4)] $F_{2\ell,d}(a) \leq F_{2\ell,d}(b)$ whenever $a\le b$.
\end{itemize}
We will see that $\mathcal{B}_d\neq \emptyset$ for every $d\ge3$ in the next section.
For a finite subset $X\subset \mathbb{S}^{d-1}$ with $A(X)\setminus\{-1\}\ne\emptyset$, we also let
\[
 a_X:=\max\bigl\{\,|a| \;\bigm|\; a\in A(X)\setminus\{-1\}\,\bigr\}.
\]
Since $\pm1\notin A(X)\setminus\{-1\}$, we have $a_X\in[0,1)$.

\begin{theorem}\label{thm:cutoff}
Let $X\subset\mathbb{S}^{d-1}$ be a finite antipodal subset with $A(X)\setminus\{-1\}\neq\emptyset$ (equivalently, $|X|\ge4$), and let $\{F_{2\ell,d}\}_{\ell\in\N}\in\mathcal{B}_d$.
If there exists $\ell_0\in\N$ such that
\begin{equation}\label{eq:ell_0}
  F_{2\ell_0,d}\bigl(a_X\bigr)<\frac{2}{|X|-2},
\end{equation}
then $2\ell\notin \Tev(X)$ for all $\ell\ge\ell_0$.
Such an $\ell_0$ always exists, by condition~{\rm (B1)}.
\end{theorem}

\begin{proof}
Suppose that $\sum_{\x,\y\in X} w(\x)\,w(\y)\,G_{2\ell,d}(\langle \x,\y\rangle)=0$ for some $\ell\in\N$ and $w\in \Wt(X)$.
By \eqref{eq:gegen_psi}, the evenness of $G_{2\ell,d}$, and the positivity of $\Phi_{X,w}$, we have
\begin{align*}
\Psi_{X,w}
&=\Bigl|\sum_{\x,\y\in X}w(\x)\,w(\y)\,G_{2\ell,d}(\langle \x,\y\rangle)-\Psi_{X,w}\Bigr| \\
&\le \sum_{a\in A(X)\setminus\{-1\}}\Phi_{X,w}(a)\,\bigl|G_{2\ell,d}(|a|)\bigr| .
\end{align*}
By conditions (B2) and (B4), we get 
\[ \bigl|G_{2\ell,d}(|a|)\bigr|\le  \max_{b\in A(X)\setminus\{-1\}} F_{2\ell,d}(|b|)=F_{2\ell,d}\bigl(a_X\bigr).\]
Thus, by Lemma \ref{lem:mass}, one has
\[ \Psi_{X,w}\le  \sum_{a\in A(X)\setminus\{-1\}}\Phi_{X,w}(a) F_{2\ell,d}\bigl(a_X\bigr) =(1-\Psi_{X,w})  F_{2\ell,d}\bigl(a_X\bigr) \]
and the factor $1-\Psi_{X,w}$ is strictly positive because $A(X)\setminus\{-1\}\ne\emptyset$.
Hence
\begin{equation}\label{eq:bound}
  \frac{\Psi_{X,w}}{1-\Psi_{X,w}}\;\le\; F_{2\ell,d}\bigl(a_X\bigr).
\end{equation}
The function $t\mapsto t/(1-t)$ is increasing on $[0,1)$, and $\Psi_{X,w}\ge\Psi_X^{\inf}=2/|X|$ by Lemma \ref{lem:Psi_inf}; therefore \eqref{eq:bound} yields
\begin{equation}\label{eq:ell_0_gen}
 \frac{2}{|X|-2}=\frac{2/|X|}{1-2/|X|}
 =\frac{\Psi_X^{\inf}}{1-\Psi_X^{\inf}}
 \;\le\; F_{2\ell,d}\bigl(a_X\bigr).
\end{equation}
Consequently, if $\ell$ is such that $F_{2\ell,d}\bigl(a_X\bigr)<2/(|X|-2)$, then $2\ell\notin \Tev(X)$.
Assuming \eqref{eq:ell_0}, condition (B3) gives $F_{2\ell,d}\bigl(a_X\bigr)\le F_{2\ell_0,d}\bigl(a_X\bigr)<2/(|X|-2)$ for every $\ell\ge\ell_0$, which completes the proof.
\end{proof}

\begin{remark}\label{rem:non_antipodal}
The assumption that $X$ be antipodal enters the proof of Theorem \ref{thm:cutoff} only through the value of $\Psi_X^{\inf}$: for an arbitrary finite subset $X\subset\mathbb{S}^{d-1}$ with $A(X)\setminus\{-1\}\ne\emptyset$, the same argument shows that $F_{2\ell_0,d}\bigl(a_X\bigr)<\Psi_X^{\inf}/(1-\Psi_X^{\inf})$ implies $2\ell\notin\Tev(X)$ for all $\ell\ge\ell_0$, and $\Psi_X^{\inf}$ can be computed in general by the argument of Lemma \ref{lem:Psi_inf}.
\end{remark}

\section{Laplace cutoff}\label{sec:effective}

For each $d\ge3$, we now show that the family $\{F^{\mathrm{Lap}}_{2\ell,d}\}_{\ell\in\N}$ defined in \eqref{eq:intro-FLap} belongs to $\mathcal{B}_d$.
Recall
\[
 s_d :=\int_0^\pi\sin^{d-3}\theta\,d\theta.
\]

\begin{lemma}
\label{lem:Lap}
Let $d\ge3$, and let $F^{\mathrm{Lap}}_{2\ell,d}$ $(\ell\in\N)$ be the functions on $[0,1)$ defined by \eqref{eq:intro-FLap}.
Then $\{F^{\mathrm{Lap}}_{2\ell,d}\}_{\ell\in\N}\in\mathcal{B}_d$.
\end{lemma}

\begin{proof}
The Laplace integral formula for the normalized Gegenbauer polynomial (see \cite[Eq.~(4.1.2)]{Askey}) gives
\[
 G_{2\ell,d}(a)=
 \frac{1}{s_d} \int_0^\pi
 \bigl(a+i\sqrt{1-a^2}\cos\theta\bigr)^{2\ell}
 \sin^{d-3}\theta\,d\theta.
\]
Since
\[
 \left|a+i\sqrt{1-a^2}\cos\theta\right|^{2\ell}
 =\bigl(a^2+(1-a^2)\cos^2\theta\bigr)^\ell,
\]
the triangle inequality proves condition~(B2).
We next prove condition~(B1).
Fix $a\in[0,1)$.
Since
\begin{equation}\label{eq:base}
 a^2+(1-a^2)\cos^2\theta
 =1-(1-a^2)\sin^2\theta,
\end{equation}
this quantity belongs to $[0,1]$.
Let $\delta\in(0,\pi/2)$.
For $\theta\in[\delta,\pi-\delta]$, we have
\[
 a^2+(1-a^2)\cos^2\theta
 \leq 1-(1-a^2)\sin^2\delta<1.
\]
On the other hand, the same quantity is at most $1$ on the two remaining intervals.
Splitting the integral into the intervals $[0,\delta]$, $[\delta,\pi-\delta]$, and $[\pi-\delta,\pi]$, we therefore obtain
\[
 F^{\mathrm{Lap}}_{2\ell,d}(a)
 \leq
 \frac{2}{s_d}\int_0^\delta
 \sin^{d-3}\theta\,d\theta
 +
 \bigl(1-(1-a^2)\sin^2\delta\bigr)^\ell.
\]
Given $\varepsilon>0$, we first choose $\delta>0$ sufficiently small so that the first term on the right-hand side is less than $\varepsilon/2$.
With this $\delta$ fixed, the number $1-(1-a^2)\sin^2\delta$ is strictly smaller than $1$.
Hence, for all sufficiently large $\ell$, the second term is also less than $\varepsilon/2$.
It follows that
\[
 \lim_{\ell\to\infty}F^{\mathrm{Lap}}_{2\ell,d}(a)=0,
\]
which proves condition~(B1).
Next, since
\[
 0\leq a^2+(1-a^2)\cos^2\theta\leq1,
\]
we have
\[
 \bigl(a^2+(1-a^2)\cos^2\theta\bigr)^{\ell+1}
 \leq
 \bigl(a^2+(1-a^2)\cos^2\theta\bigr)^\ell.
\]
Integrating this inequality proves condition~(B3).
Finally, by \eqref{eq:base}, for each fixed $\theta$ the quantity $a^2+(1-a^2)\cos^2\theta=1-(1-a^2)\sin^2\theta$ is nonnegative and nondecreasing in $a\in[0,1)$; hence so is its $\ell$-th power.
Integrating in $\theta$ proves condition~(B4).
\end{proof}

\begin{remark}\label{rem:F_Lap_closed}
Using
\[
 \frac{1}{s_d}\int_0^\pi\cos^{2k}\theta\,\sin^{d-3}\theta\,d\theta
 =\frac{\bigl(\tfrac12\bigr)_k}{\bigl(\tfrac{d-1}{2}\bigr)_k},
\]
where $(x)_k:=x(x+1)\cdots(x+k-1)$ denotes the Pochhammer symbol, we obtain
\begin{equation}\label{eq:F_Lap_closed}
 F^{\mathrm{Lap}}_{2\ell,d}(a)
 =\sum_{k=0}^{\ell}\binom{\ell}{k}
 \frac{\bigl(\tfrac12\bigr)_k}{\bigl(\tfrac{d-1}{2}\bigr)_k}\,
 a^{2(\ell-k)}\bigl(1-a^2\bigr)^{k},
\end{equation}
a finite sum with nonnegative terms.
Note that the inequality $|G_{2\ell,d}(a)|\le F^{\mathrm{Lap}}_{2\ell,d}(a)$ in condition~(B2) is an equality at $a=0$, namely $F^{\mathrm{Lap}}_{2\ell,d}(0) =\bigl(\tfrac12\bigr)_\ell\big/\bigl(\tfrac{d-1}{2}\bigr)_\ell =|G_{2\ell,d}(0)|$ for all $\ell\in\N$.
\end{remark}

We are now in a position to prove Theorem \ref{thm:main}.
\begin{proof}[Proof of Theorem \ref{thm:main}]
If $|X|\le2$, then $A(X)\setminus\{-1\}=\emptyset$, so that \eqref{eq:gegen_psi} and Lemma \ref{lem:mass} give
\[\sum_{\x,\y\in X}w(\x)w(\y)G_{2\ell,d}(\langle\x,\y\rangle)=\Psi_{X,w}=1\ne0\]
for every $w\in\Wt(X)$; hence $\Tev(X)=\emptyset$.
If $|X|\ge4$, then $A(X)\setminus\{-1\}\ne\emptyset$, and by Lemma \ref{lem:Lap} we have $\{F^{\mathrm{Lap}}_{2\ell,d}\}_{\ell\in\N}\in\mathcal{B}_d$, so that $\ell_0$ exists by condition (B1).
Theorem \ref{thm:cutoff} then gives the assertion.
\end{proof}

Sharper pointwise bounds for Gegenbauer polynomials, such as those of Krasikov~\cite{Kra17}, may give a smaller value of $\ell_0$ in particular examples.
Lemma~\ref{lem:Lap}, however, has the advantage of giving a simple family in $\mathcal{B}_d$ for every $d\geq3$ without any additional restriction on $a$.

\begin{remark}

In the equal-weight setting, Misawa and Nishimura~\cite{MN26} gave an effective proof that a finite set $X\subset\mathbb{S}^{d-1}$ has infinite harmonic strength if and only if it is antipodal, and obtained an explicit $\ell_1$ such that $2\ell\notin T(X)$ for all $\ell\ge\ell_1$.
Throughout this remark we have translated their statements into the notation of the present paper; their $d$-dimensional sphere $\mathbb{S}^{d}$ is our $\mathbb{S}^{d-1}$.
For antipodal $X$, their argument in fact yields the same bound for $\Tev(X)$.
Indeed, their pointwise Gegenbauer estimate \cite[Lemma~2.5 and Theorem~2.7]{MN26} applies to the weighted sum in \eqref{eq:gegen_psi}, and Lemmas~\ref{lem:mass} and~\ref{lem:Psi_inf} give, for $|X|\ge4$,
\[
 \frac{1-\Psi_{X,w}}{\Psi_{X,w}}
 \le\frac{|X|-2}{2}
 \qquad\text{for every }w\in\Wt(X).
\]
The right-hand side is exactly the corresponding ratio for equal weights in their argument.
Their bound on $\ell$ is therefore uniform over positive weights, even when the weight depends on the index.
The family $\{F^{\mathrm{Lap}}_{2\ell,d}\}_{\ell\in\N}$ above provides an alternative explicit bound, computable from the finite sum~\eqref{eq:F_Lap_closed}.
The results of the present paper and those of~\cite{MN26} were obtained independently.
\end{remark}

\section{Elimination}\label{sec:elimination}
For an antipodal subset $X\subset\mathbb{S}^{d-1}$, once the Laplace cutoff $\ell_0$ of $X$ has been computed, there remains the problem of deciding, for each of the finitely many even integers $2\ell\in\{2,4,\ldots,2\ell_0-2\}$, whether $2\ell\in\Tev(X)$.
In this section, we present an approach to this problem for finite subsets preserved by a finite subgroup of the orthogonal group.
This setting covers unions of shells of a lattice, since the automorphism group of a lattice preserves each of its shells.

Although we assume $d\ge3$ throughout the paper, the results of this section hold for $d\ge2$.
Throughout this section, let $X\subset\mathbb{S}^{d-1}$ be a finite subset, not necessarily antipodal, and let $H\subset O(\R^d)$ be a finite subgroup preserving $X$.
Let
\[ X=\mathcal{O}_1\sqcup\cdots\sqcup\mathcal{O}_r \]
be the decomposition of $X$ into $H$-orbits.
For each $\ell\in\N$, we now define an $r\times r$ matrix $M_\ell(X)$; although it depends on $H$, we suppress $H$ from the notation.
Let $\varphi_{\ell,1},\ldots,\varphi_{\ell,h_\ell}$ be an orthonormal basis of $\mathrm{Harm}_\ell$ with respect to the inner product \eqref{eq:harm-inner}.
For each $H$-orbit $\mathcal{O}_i$, define the vector
\[ \boldsymbol{v}_\ell(\mathcal{O}_i):=\left(\frac{1}{|\mathcal{O}_i|}\sum_{\x\in\mathcal{O}_i}\varphi_{\ell,\nu}(\x)\right)_{1\le\nu\le h_\ell}\in\R^{h_\ell},\qquad i=1,2,\ldots,r, \]
and let $M_\ell(X)$ be the Gram matrix of these vectors:
\[ M_\ell(X):=\bigl(\langle\boldsymbol{v}_\ell(\mathcal{O}_i),\boldsymbol{v}_\ell(\mathcal{O}_j)\rangle\bigr)_{1\le i,j\le r}. \]
Being a Gram matrix, $M_\ell(X)$ is symmetric and positive semidefinite.
By the addition formula \eqref{eq:addition}, the $(i,j)$-entry of $M_\ell(X)$ can be written as
\begin{equation}\label{eq:v_inner}
 \langle\boldsymbol{v}_\ell(\mathcal{O}_i),\boldsymbol{v}_\ell(\mathcal{O}_j)\rangle
 =\frac{h_\ell}{|\mathcal{O}_i||\mathcal{O}_j|}\sum_{\x\in\mathcal{O}_i,\,\y\in\mathcal{O}_j}G_{\ell,d}(\langle\x,\y\rangle).
\end{equation}
In particular, $M_\ell(X)$ does not depend on the choice of the orthonormal basis.
For $\boldsymbol{c}=(c_1,\ldots,c_r)\in\R^r$, we also consider the quadratic form $M_\ell(X)[\boldsymbol{c}]:=\boldsymbol{c}\,M_\ell(X)\,{}^t\boldsymbol{c}$.
By definition, we have
\begin{equation}\label{eq:quad}
 M_\ell(X)[\boldsymbol{c}]
 =\left\langle\sum_{i=1}^r c_i\boldsymbol{v}_\ell(\mathcal{O}_i),\ \sum_{i=1}^r c_i\boldsymbol{v}_\ell(\mathcal{O}_i)\right\rangle
 =\left\|\sum_{i=1}^r c_i\boldsymbol{v}_\ell(\mathcal{O}_i)\right\|^2\ge0.
\end{equation}

Together with Proposition~\ref{prop:averaging}, the following proposition characterizes the condition $\ell\in T_\ast(X)$ in terms of the matrix $M_\ell(X)$.

\begin{proposition}\label{prop:orbit_matrix}
Let $X\subset\mathbb{S}^{d-1}$ be a finite subset, let $H\subset O(\R^d)$ be a finite subgroup preserving $X$, and let $X=\mathcal{O}_1\sqcup\cdots\sqcup\mathcal{O}_r$ be the decomposition of $X$ into $H$-orbits.
Then, for $\ell\in\N$ and $\boldsymbol{c}=(c_1,\ldots,c_r)\in\R^r$, the following conditions are equivalent{\rm :}
\begin{itemize}
\item[$(\mathrm{i})$] $\displaystyle\sum_{i=1}^r\frac{c_i}{|\mathcal{O}_i|}\sum_{\x\in\mathcal{O}_i}P(\x)=0$ for all $P\in\mathrm{Harm}_\ell${\rm ;}
\item[$(\mathrm{ii})$] $M_\ell(X)[\boldsymbol{c}]=0$, that is, the quadratic form associated with $M_\ell(X)$ vanishes at $\boldsymbol{c}${\rm ;}
\item[$(\mathrm{iii})$] $\boldsymbol{c}\,M_\ell(X)=\boldsymbol{0}$, that is, $\boldsymbol{c}$ lies in the kernel of $M_\ell(X)$.
\end{itemize}
\end{proposition}
\begin{proof}
Put $\boldsymbol{u}:=\sum_{i=1}^r c_i\boldsymbol{v}_\ell(\mathcal{O}_i)\in\R^{h_\ell}$.
By \eqref{eq:quad}, condition $(\mathrm{ii})$ is equivalent to $\boldsymbol{u}=\boldsymbol{0}$.
Since the $j$-th entry of $\boldsymbol{c}\,M_\ell(X)$ equals $\langle\boldsymbol{u},\boldsymbol{v}_\ell(\mathcal{O}_j)\rangle$, the equality $\boldsymbol{u}=\boldsymbol{0}$ implies $(\mathrm{iii})$; conversely, $(\mathrm{iii})$ implies $M_\ell(X)[\boldsymbol{c}]=\boldsymbol{c}\,M_\ell(X)\,{}^t\boldsymbol{c}=0$, that is, $(\mathrm{ii})$.
Finally, the $\nu$-th entry of $\boldsymbol{u}$ is the left-hand side of $(\mathrm{i})$ with $P=\varphi_{\ell,\nu}$.
Since $\varphi_{\ell,1},\ldots,\varphi_{\ell,h_\ell}$ form a basis of $\mathrm{Harm}_\ell$ and the left-hand side of $(\mathrm{i})$ is linear in $P$, condition $(\mathrm{i})$ is also equivalent to $\boldsymbol{u}=\boldsymbol{0}$.
\end{proof}

Combining Proposition~\ref{prop:orbit_matrix} with Proposition~\ref{prop:averaging}, we obtain the following characterization of $T_\ast(X)$ in terms of the matrices $M_\ell(X)$.

\begin{corollary}\label{cor:kernel_criterion}
Let $X$, $H$ and $\mathcal{O}_1,\ldots,\mathcal{O}_r$ be as in Proposition~\ref{prop:orbit_matrix}, and let $\ell\in\N$.
\begin{itemize}
\item[$(\mathrm{i})$] $\ell\in T_\ast(X)$ if and only if $\ker M_\ell(X)$ contains a vector with all entries positive.
In this case, for any such vector $\boldsymbol{c}=(c_1,\ldots,c_r)$ normalized so that $c_1+\cdots+c_r=1$, the weight $w$ defined by $w(\x):=c_i/|\mathcal{O}_i|$ for $\x\in\mathcal{O}_i$ realizes $\ell$.
\item[$(\mathrm{ii})$] In particular, if $\det M_\ell(X)\ne0$, then $\ell\notin T_\ast(X)$.
\end{itemize}
\end{corollary}
\begin{proof}
(i) Suppose that $\ell\in T_\ast(X)$.
By Proposition~\ref{prop:averaging}, there exists a vector $\boldsymbol{c}=(c_1,\ldots,c_r)\in \R_{>0}^r$ with $c_1+\cdots+c_r=1$ satisfying $\sum_{i=1}^r\frac{c_i}{|\mathcal{O}_i|}\sum_{\x\in\mathcal{O}_i}P(\x)=0$ for all $P\in\mathrm{Harm}_\ell$. 
Thus, $\boldsymbol{c}\in \ker M_\ell(X)$ because of Proposition~\ref{prop:orbit_matrix}.
Conversely, let $\boldsymbol{c}\in\ker M_\ell(X)$ have all entries positive.
Since $\ker M_\ell(X)$ is a linear subspace, we may assume $c_1+\cdots+c_r=1$ after multiplying $\boldsymbol{c}$ by a positive constant.
By the equivalence of (i) and (iii) in Proposition~\ref{prop:orbit_matrix}, the numbers $c_1,\ldots,c_r$ satisfy \eqref{eq:orbit_condition} for every $P\in\mathrm{Harm}_\ell$, and hence $\ell\in T_\ast(X)$ by Proposition~\ref{prop:averaging}.

(ii) If $\det M_\ell(X)\ne0$, then $\ker M_\ell(X)=\{\boldsymbol{0}\}$, which contains no vector with all entries positive.
Hence $\ell\notin T_\ast(X)$ by (i).
\end{proof}

\begin{remark}\label{rem:h-orbit}
In the same vein, one can show that $\ell\in T(X)$ if and only if the vector $(|\mathcal{O}_1|,\ldots,|\mathcal{O}_r|)$ lies in $\ker M_\ell(X)$.
\end{remark}

In order to apply Corollary~\ref{cor:kernel_criterion}, one has to compute the $r\times r$ matrix $M_\ell(X)$ explicitly.
We now describe how to compute it for a given $\ell$.

For $i,j\in\{1,2,\ldots,r\}$, define
\[ A(\mathcal{O}_i,\mathcal{O}_j):=\{\langle\x,\y\rangle\mid\x\in\mathcal{O}_i,\ \y\in\mathcal{O}_j\}, \]
and, for $a\in A(\mathcal{O}_i,\mathcal{O}_j)$, set
\[ n_{ij}(a):=\bigl|\{(\x,\y)\in\mathcal{O}_i\times\mathcal{O}_j\mid\langle\x,\y\rangle=a\}\bigr|. \]
Note that, unlike $A(X)$, the set $A(\mathcal{O}_i,\mathcal{O}_i)$ contains $1$, since pairs with $\x=\y$ are not excluded.
We call the collection $\{(a,n_{ij}(a))\mid a\in A(\mathcal{O}_i,\mathcal{O}_j)\}$ the \emph{inner product distribution between $\mathcal{O}_i$ and $\mathcal{O}_j$}.
Once the inner product distributions have been computed, by \eqref{eq:v_inner}, the $(i,j)$-entry of $M_\ell(X)$ is given by
\begin{equation}\label{eq:dist_dist}
 \langle\boldsymbol{v}_\ell(\mathcal{O}_i),\boldsymbol{v}_\ell(\mathcal{O}_j)\rangle= \frac{h_\ell}{|\mathcal{O}_i||\mathcal{O}_j|} \sum_{a\in A(\mathcal{O}_i,\mathcal{O}_j)}n_{ij}(a)\,G_{\ell,d}(a).
\end{equation}
Since the inner product distributions do not depend on $\ell$, they need to be computed only once for all $\ell$.

\section{The cubic lattice $\Z^3$}\label{sec:Z3}

Let $\Lambda=\Z^3$.
We first recall the geometry and harmonic strengths of its normalized shells $\Lambda[1]$ and $\Lambda[3]$.
These are the vertex sets of a regular octahedron and a cube.
With $\boldsymbol{e}_1,\boldsymbol{e}_2,\boldsymbol{e}_3$ denoting the standard basis vectors of $\R^3$, set
\begin{align*}
 \mathcal{O}_1&:=\Lambda[1]=\{\pm \boldsymbol{e}_1,\pm \boldsymbol{e}_2,\pm \boldsymbol{e}_3\},\\
 \mathcal{O}_2&:=\Lambda[3] =\frac1{\sqrt3}\{(\varepsilon_1,\varepsilon_2,\varepsilon_3) \mid\varepsilon_i\in\{\pm1\}\}.
\end{align*}
The orthogonal automorphism group of $\Z^3$ is the Weyl group $W(B_3)\cong (\Z/2\Z)^3\rtimes S_3$ of order $48$, which consists of the signed permutations of the coordinates.
The shells $\mathcal{O}_1$ and $\mathcal{O}_2$ are single $W(B_3)$-orbits:
\[
 \mathcal{O}_1=W(B_3)\boldsymbol{e}_1,\qquad \mathcal{O}_2=W(B_3)\boldsymbol{u}
\]
with $\boldsymbol{u}=\frac{1}{\sqrt{3}}(\boldsymbol{e}_1+\boldsymbol{e}_2+\boldsymbol{e}_3)$.
They are distinct orbits because the number of nonzero coordinates is preserved by $W(B_3)$.

Both $\mathcal{O}_1$ and $\mathcal{O}_2$ are classical spherical $3$-designs, and their ordinary design strength is exactly $3$ (cf.~\cite{DGS77}); see also Hardin--Sloane~\cite[Table~I]{HS96}.
Their full harmonic strengths are $T(\mathcal{O}_1)=T(\mathcal{O}_2)=\{2\}\cup\mathbb{O}$.
Indeed, since $|\mathcal{O}_1|=6$, $a_{\mathcal{O}_1}=0$, $|\mathcal{O}_2|=8$ and $a_{\mathcal{O}_2}=\frac13$, the Laplace cutoffs of $\mathcal{O}_1$ and $\mathcal{O}_2$ are $2$ and $4$, respectively, so that by Theorem~\ref{thm:main} only $k=2$ for $\mathcal{O}_1$ and $k\in\{2,4,6\}$ for $\mathcal{O}_2$ remain to be examined.
These cases are easily checked by direct computation.

We now consider their union, for which positive weights can enlarge the harmonic strength.
The following theorem determines the complete set $T_\ast(\mathcal{O}_1\sqcup \mathcal{O}_2)$, including the additional indices $4$ and $10$.

\begin{theorem}\label{thm:Z3-shells}
Let $\Lambda=\Z^3$, let $X=\Lambda[1,3]$, and let $\ell_0$ be the Laplace cutoff of $X$.
Then
\[
 \ell_0=18,\qquad
 T(X)=\{2\}\cup\mathbb{O},\qquad
 T_\ast(X)=\{2,4,10\}\cup\mathbb{O}.
\]
Moreover, $X$ is not a weighted $\{4,10\}$-design: these two indices require different weight functions.
\end{theorem}

\begin{proof}
Recall that $X=\mathcal{O}_1\sqcup\mathcal{O}_2$ with $\mathcal{O}_1=\Lambda[1]$ and $\mathcal{O}_2=\Lambda[3]$.
We have $|\mathcal{O}_1|=6$ and $|\mathcal{O}_2|=8$, so that $|X|=14$.
One can compute
\[ A(X)=\{-1,0,\pm\frac13,\pm\frac{1}{\sqrt3}\},\]
and hence, $a_X=1/\sqrt3$.
Exact evaluation of~\eqref{eq:F_Lap_closed} gives
\begin{align*}
 F^{\mathrm{Lap}}_{36,3}\left(\frac{1}{\sqrt3}\right)
 =\frac{1378999684337}{8463329722368}<\frac16
 =\frac2{|X|-2}<  F^{\mathrm{Lap}}_{34,3}\left(\frac{1}{\sqrt3}\right)
 =\frac{13139436329}{78364164096}.
\end{align*}
Hence $\ell_0=18$ and $\Tev(X)\subset\{2,4,\ldots,34\}$.

We now apply Corollary~\ref{cor:kernel_criterion} with $H:=W(B_3)$, which preserves $X$ and has exactly two orbits $\mathcal{O}_1$ and $\mathcal{O}_2$ on $X$.
The inner product distributions are given as follows:
\[
\begin{array}{c|l}
 (i,j) & \{(a,n_{ij}(a))\mid a\in A(\mathcal{O}_i,\mathcal{O}_j)\}\\ \hline
 (1,1) & (1,6),\ (0,24),\ (-1,6)\\
 (2,2) & (1,8),\ (\tfrac13,24),\ (-\tfrac13,24),\ (-1,8)\\
 (1,2) & (\tfrac1{\sqrt3},24),\ (-\tfrac1{\sqrt3},24)
\end{array}
\]

Let $k$ be an even positive integer.
Since $G_{k,3}$ is an even polynomial and $h_k=\dim\mathrm{Harm}_k=2k+1$ for $d=3$, it follows from $G_{k,3}(\pm1)=1$ and \eqref{eq:dist_dist} that
\begin{equation}\label{eq:Z3-Mk}
 M_k(X)=(2k+1)\begin{pmatrix}\alpha_k&\beta_k\\\beta_k&\gamma_k\end{pmatrix}
\end{equation}
with
\[
 \alpha_k=\frac{1+2G_{k,3}(0)}{3},\qquad
\beta_k=G_{k,3}\Bigl(\frac1{\sqrt3}\Bigr),\qquad
 \gamma_k=\frac{1+3G_{k,3}(\frac13)}{4}.
\]
All these entries are rational and can, for instance, be evaluated using
\[
 G_{2m,3}(\sqrt r)=\frac1{2^{2m}}\sum_{j=0}^m
 (-1)^j\frac{(4m-2j)!}{j!(2m-j)!(2m-2j)!}\,r^{m-j}
 \qquad (r\ge0).
\]
Exact calculation shows that $\det M_k(X)>0$ for $k=12,16,18,20,\ldots,34$, so that none of these $k$ belongs to $\Tev(X)$ by Corollary~\ref{cor:kernel_criterion}~(ii).
For the remaining even integers $k=2,4,6,8,10,14$, we have $\det M_k(X)=0$, and
\[
\begin{array}{c|cccccc}
 k&2&4&6&8&10&14\\ \hline
 \alpha_k&0&\frac7{12}&\frac18&\frac{33}{64}&\frac{65}{384}&\frac{595}{3072}\\[3pt]
 \beta_k&0&-\frac7{18}&\frac29&\frac{11}{72}&-\frac{65}{243}&\frac{595}{2916}\\[3pt]
 \gamma_k&0&\frac7{27}&\frac{32}{81}&\frac{11}{243}&\frac{8320}{19683}&\frac{38080}{177147}
\end{array}
\]
For $k=2$, the matrix $M_2(X)$ is zero, so that $\ker M_2(X)=\R^2$ and $2\in\Tev(X)$.
For $k\in\{4,6,8,10,14\}$, we always have $\alpha_k>0$, so that $M_k(X)$ has rank one and
\[
 \ker M_k(X)=\R\,(-\beta_k,\alpha_k).
\]
The space $ \ker M_k(X)$ contains a vector with both entries positive if and only if $\beta_k<0$, which is the case exactly for $k=4$ and $k=10$.
Therefore, Corollary~\ref{cor:kernel_criterion}~(i) gives $\Tev(X)=\{2,4,10\}$.

By Remark \ref{rem:h-orbit}, an even integer $k$ belongs to $T(X)$ if and only if $(6,8)\in\ker M_k(X)$.
Since $(6,8)$ has positive entries, this is possible only for $k\in\{2,4,10\}$, and $(6,8)$ is not proportional to $(-\beta_k,\alpha_k)$ for $k=4,10$.
Hence $T(X)\cap2\N=\{2\}$, and $T(X)=\{2\}\cup\mathbb{O}$.

Finally, the last statement follows from $\ker M_4(X)\cap\ker M_{10}(X)=\{\boldsymbol{0}\}$.
\end{proof}

\begin{remark}\label{rem:Z3-weights}
For $\boldsymbol{c}=(c_1,c_2)\in\R_{>0}^2$ with $c_1+c_2=1$, let $w_{\boldsymbol{c}}\in\Wt(X)$ be the weight defined by $w_{\boldsymbol{c}}(\x):=c_1/6$ for $\x\in\mathcal{O}_1$ and $w_{\boldsymbol{c}}(\x):=c_2/8$ for $\x\in\mathcal{O}_2$.
By the proof of Theorem~\ref{thm:Z3-shells}, the weights of this form realizing the even indices of $X$ are as follows:
\[
\begin{array}{c|c|c|c}
 k & \boldsymbol{c}=(c_1,c_2) & w_{\boldsymbol{c}}(\x),\ \x\in\mathcal{O}_1 & w_{\boldsymbol{c}}(\x),\ \x\in\mathcal{O}_2\\ \hline
 2 & \text{arbitrary} & c_1/6 & c_2/8\\[2pt]
 4 & \bigl(\frac25,\frac35\bigr) & \frac1{15} & \frac3{40}\\[4pt]
 10 & \bigl(\frac{128}{209},\frac{81}{209}\bigr) & \frac{64}{627} & \frac{81}{1672}
\end{array}
\]
For $k=2$, every $w_{\boldsymbol{c}}$ realizes $k$ since $M_2(X)=O$; this includes the equal weight $w\equiv1/14$, which corresponds to $\boldsymbol{c}=(\frac37,\frac47)$.
For $k=4$ and $k=10$, the vector $\boldsymbol{c}$ is the normalization of $(-\beta_k,\alpha_k)$, which spans $\ker M_k(X)$.
\end{remark}

We remark that the weight $w_{\boldsymbol{c}}$ with $\boldsymbol{c}=(\frac25,\frac35)$ realizes every $\ell\le5$: the even ones by Remark~\ref{rem:Z3-weights}, and the odd ones since $w_{\boldsymbol{c}}(-\x)=w_{\boldsymbol{c}}(\x)$.
Hence $X$ is a weighted $5$-design; this recovers the weighted $5$-design on $\Lambda[1,3]$ obtained from Bajnok's construction~\cite[Corollary~9(1) and the subsequent remark]{Baj07}.

\appendix
\section{The cubature formulas of de la Harpe, Pache and Venkov}\label{app:HPV}

In this appendix we recall the results of de la Harpe, Pache and Venkov~\cite[Theorems~$D_4$ and~$E_8$]{HPV07} in dimensions $4$ and $8$, and determine the harmonic strength and the weighted harmonic strength of the sets they treated.
In both cases the subscript $r$ of the sets $X_r$ and $Y_r$ below is the number of their orbits under the relevant group $H$.

\subsection{Shells of $D_4\cup D_4'$}\label{app:D4}
Let
\[
 D_4:=\{\x\in\Z^4\mid x_1+x_2+x_3+x_4\in2\Z\}
\]
be the root lattice of type $D_4$, let $D_4^\ast:=\{\y\in\R^4\mid\langle\x,\y\rangle\in\Z\ \text{for all}\ \x\in D_4\}$ be its dual lattice, and put
\[
 D_4':=\sqrt2\,D_4^\ast,\qquad \Lambda:=D_4\cup D_4'.
\]
The lattice $D_4'$ is isometric to $D_4$, and $D_4\cap D_4'=\{\boldsymbol{0}\}$, since the nonzero coordinates of the vectors of $D_4'$ are irrational.
Although $\Lambda$ is not a lattice, we use the notation~\eqref{eq:shell} for its normalized shells.
The numbers of vectors of $D_4$ of squared norms $2,4,6,8,10$ are $24,24,96,24,144$, and the same holds for $D_4'$.
Hence, we obtain
\[
 |\Lambda[2]|=48,\qquad |\Lambda[6]|=192,\qquad |\Lambda[10]|=288.
\]
Since $\sqrt2\,D_4\subset D_4'$ and $\sqrt2\,D_4'=2D_4^\ast\subset D_4$, we have $\Lambda[n]\subset\Lambda[2n]$ for every $n$; for instance $\Lambda[4]=\Lambda[2]$, since $\Lambda[2]\subset\Lambda[4]$ and both sets consist of $48$ points (see~\cite[Section~6.4]{HPV07}).
In particular a union such as $\Lambda[2,4]$ gives nothing new, and the constructions below use the shells of squared norms $2$, $6$ and $10$.

Let $H$ be the group of orthogonal transformations of $\R^4$ preserving $\Lambda$.
As described in~\cite[Section~6.4]{HPV07}, $H$ is generated by the automorphism group $\operatorname{Aut}(D_4)=W(F_4)$ of $D_4$ and an orthogonal transformation interchanging $D_4$ and $D_4'$, for instance
\[
 U:=\frac1{\sqrt2}
 \begin{pmatrix}
  1&1&0&0\\-1&1&0&0\\0&0&1&1\\0&0&-1&1
 \end{pmatrix}.
\]
The group $H$ preserves every normalized shell of $\Lambda$, and each of the sets
\[
 \mathcal{O}_1:=\Lambda[2],\qquad
 \mathcal{O}_2:=\Lambda[6],\qquad
 \mathcal{O}_3:=\Lambda[10]
\]
is a single $H$-orbit; these three orbits are pairwise disjoint.

\begin{theorem}[{\cite[Theorem~$D_4$]{HPV07}}]\label{thm:HPV-D4}
\begin{itemize}
\item[$(\mathrm{i})$] For every $m\in\N$, the set $\Lambda[2m]$ is a spherical $7$-design.
In particular, $X_1:=\Lambda[2]$ is a spherical $7$-design with $|X_1|=48$.
\item[$(\mathrm{ii})$] $X_2:=\Lambda[2,6]$ is a weighted $11$-design with $|X_2|=240$.
\item[$(\mathrm{iii})$] $X_3:=\Lambda[2,6,10]$ is a weighted $15$-design with $|X_3|=528$.
\end{itemize}
\end{theorem}

Assertion~(i) for $m=1$ goes back to Goethals--Seidel~\cite[Section~6]{GS81}.
Under the standard identification of $\R^4$ with the quaternions, $X_1$ is the binary octahedral group of order $48$.
Its harmonic strength is therefore covered by~\cite{HNT}.

Let us state the complete description of both $T(X_r)$ and $T_\ast(X_r)$ for $r=1,2,3$.

\begin{theorem}\label{thm:D4-shells}
Let $X_1$, $X_2$ and $X_3$ be as in Theorem~\ref{thm:HPV-D4}, and put
\[
 S_D:=\{2,4,6,10,14,22\}.
\]
Then the Laplace cutoffs of $X_1$, $X_2$ and $X_3$ are $24$, $540$ and $6600$, respectively,
\[
 T(X_1)=T(X_2)=T(X_3)=S_D\cup\mathbb{O},
\]
and
\begin{align*}
\Tev(X_1)&=S_D,\\
\Tev(X_2)&=\Tev(X_1)\cup\{8,18,26,30,38,46\},\\
\Tev(X_3)&=\Tev(X_2)\cup\{12,34,42,54,62,70\}.
\end{align*}
\end{theorem}

For $X_1$, the equality $T(X_1)=T_\ast(X_1)$ is a consequence of Proposition~\ref{prop:averaging}, since $X_1$ is a single $H$-orbit.

\begin{proof}[Outline of the proof of Theorem~\ref{thm:D4-shells}]
By~\eqref{eq:v_inner}, the $(i,j)$-entry of $M_k(X_r)$ depends only on $\mathcal{O}_i$ and $\mathcal{O}_j$.
Hence, with the orbits ordered as $\mathcal{O}_1,\mathcal{O}_2,\mathcal{O}_3$, the matrices $M_k(X_1)$ and $M_k(X_2)$ are the leading principal submatrices of $M_k(X_3)$ of sizes $1$ and $2$, respectively.

\emph{Step 1: inner product distributions.}
For even $k$ the polynomial $G_{k,4}$ is even, so that only the squares of the inner products matter.
For $\x\in\mathcal{O}_i$ and $b\in[0,1]$, put
\[
 m_{ij}(b):=\bigl|\{\y\in\mathcal{O}_j\mid\langle\x,\y\rangle^2=b\}\bigr|,
\]
which does not depend on the choice of $\x$, since $H$ acts transitively on $\mathcal{O}_i$ and preserves $\mathcal{O}_j$.
Then $n_{ij}(\sqrt b)+n_{ij}(-\sqrt b)=|\mathcal{O}_i|\,m_{ij}(b)$ for $b>0$ and $n_{ij}(0)=|\mathcal{O}_i|\,m_{ij}(0)$, and since $h_k=(k+1)^2$ for $d=4$, \eqref{eq:dist_dist} becomes
\begin{equation}\label{eq:D4-Mk}
 \bigl(M_k(X_3)\bigr)_{ij}=\frac{(k+1)^2}{|\mathcal{O}_j|}\sum_{b}m_{ij}(b)\,G_{k,4}(\sqrt b)
 \qquad(k\in2\N).
\end{equation}
Table~\ref{tab:D4-dist} lists the values $b$ together with $m_{ij}(b)$ for $i\le j$; the remaining ones follow from $|\mathcal{O}_i|\,m_{ij}(b)=|\mathcal{O}_j|\,m_{ji}(b)$.
Altogether $49$ distinct values of $b$ occur.

\begin{table}[ht]
\centering
\small
$\begin{array}{c|l}
 (i,j) & b\ \bigl(m_{ij}(b)\bigr)\\ \hline
 (1,1) & 0\,(18),\ \tfrac14\,(16),\ \tfrac12\,(12),\ 1\,(2)\\[2pt]
 (1,2) & 0\,(32),\ \tfrac1{12}\,(48),\ \tfrac16\,(48),\ \tfrac13\,(24),\ \tfrac23\,(24),\ \tfrac34\,(16)\\[2pt]
 (1,3) & 0\,(48),\ \tfrac1{20}\,(48),\ \tfrac1{10}\,(60),\ \tfrac15\,(12),\ \tfrac25\,(48),\ \tfrac9{20}\,(48),\ \tfrac45\,(12),\ \tfrac9{10}\,(12)\\[2pt]
 (2,2) & 0\,(30),\ \tfrac1{36}\,(24),\ \tfrac1{18}\,(24),\ \tfrac19\,(6),\ \tfrac29\,(24),\ \tfrac14\,(28),\ \tfrac49\,(12),\ \tfrac12\,(24),\ \tfrac{25}{36}\,(12),\ \tfrac89\,(6),\ 1\,(2)\\[2pt]
 (2,3) & 0\,(12),\ \tfrac1{60}\,(36),\ \tfrac1{30}\,(36),\ \tfrac1{15}\,(24),\ \tfrac2{15}\,(36),\ \tfrac3{20}\,(24),\ \tfrac4{15}\,(12),\\
       & \tfrac3{10}\,(24),\ \tfrac5{12}\,(24),\ \tfrac8{15}\,(24),\ \tfrac35\,(12),\ \tfrac{49}{60}\,(12),\ \tfrac56\,(12)\\[2pt]
 (3,3) & 0\,(34),\ \tfrac1{100}\,(24),\ \tfrac1{50}\,(18),\ \tfrac1{25}\,(8),\ \tfrac2{25}\,(16),\ \tfrac9{100}\,(16),\ \tfrac4{25}\,(16),\ \tfrac9{50}\,(32),\ \tfrac14\,(32),\\
       & \tfrac8{25}\,(16),\ \tfrac9{25}\,(2),\ \tfrac{49}{100}\,(16),\ \tfrac12\,(20),\ \tfrac{16}{25}\,(10),\ \tfrac{18}{25}\,(16),\ \tfrac{81}{100}\,(8),\ \tfrac{49}{50}\,(2),\ 1\,(2)
\end{array}$\vspace{5mm}
\caption{The squared inner products $b$ and the numbers $m_{ij}(b)$ for $\Lambda=D_4\cup D_4'$.}
\label{tab:D4-dist}
\end{table}

\emph{Step 2: the Laplace cutoffs.}
Reading off the largest $b<1$ in Table~\ref{tab:D4-dist} among the pairs $(i,j)$ with $i,j\le r$, we obtain $a_{X_1}^2=\frac12$, $a_{X_2}^2=\frac89$ and $a_{X_3}^2=\frac{49}{50}$.
For $d=4$ we have $(\frac12)_k/(\frac32)_k=1/(2k+1)$, and $f_\ell:=F^{\mathrm{Lap}}_{2\ell,4}(\sqrt b)$ satisfies
\[
 f_0=1,\qquad (2\ell+1)f_\ell=1+2\ell b\,f_{\ell-1}\quad(\ell\ge1),
\]
as one sees by substituting $u=\cos\theta$ in~\eqref{eq:intro-FLap} and integrating the derivative of $u\bigl(b+(1-b)u^2\bigr)^\ell$ over $[0,1]$.
Exact evaluation gives Table~\ref{tab:D4-cutoff} below.
\begin{table}[ht]
\centering
\small
$\begin{array}{c|ccc|ccc}
 & |X| & a_X^2 & \ell_0 & F_{2\ell_0} & \frac{2}{|X|-2} & F_{2\ell_0-2}\\ \hline
 X_1 & 48 & \tfrac12 & 24 & 41.7509\ldots & 43.4782\ldots & 43.5748\ldots\\[2pt]
 X_2 & 240 & \tfrac89 & 540 & 8.38878\ldots & 8.40336\ldots & 8.40445\ldots\\[2pt]
 X_3 & 528 & \tfrac{49}{50} & 6600 & 3.80181\ldots & 3.80228\ldots & 3.80239\ldots
\end{array}$\vspace{5mm}
\caption{The Laplace cutoffs of $X_1$, $X_2$ and $X_3$. Here $F_{2\ell}:=F^{\mathrm{Lap}}_{2\ell,4}(a_X)$, and the last three columns are multiplied by $10^3$; in each case $F_{2\ell_0}<\frac{2}{|X|-2}\le F_{2\ell_0-2}$.}
\label{tab:D4-cutoff}
\end{table}
Hence, by Theorem~\ref{thm:main}, only the even integers $k\le2\ell_0-2$, namely $k\le46$, $k\le1078$ and $k\le13198$, remain to be examined.

\emph{Step 3: evaluation of $M_k(X_3)$.}
Since $G_{k,4}(\cos\theta)=\sin((k+1)\theta)/((k+1)\sin\theta)$, for $b=p/q$ in lowest terms we have
\[
 G_{2m,4}(\sqrt b)=\frac{\rho_m(b)}{(2m+1)\,q^m},
\]
where the integers $\rho_m(b)$ are determined by
\[
 \rho_0(b)=1,\qquad \rho_1(b)=4p-q,\qquad
 \rho_{m+1}(b)=(4p-2q)\,\rho_m(b)-q^2\rho_{m-1}(b).
\]
All denominators $q$ in Table~\ref{tab:D4-dist} divide $900$.
Multiplying~\eqref{eq:D4-Mk} by $(2m+1)\,900^m/(k+1)^2$ and the $j$-th column by $|\mathcal{O}_j|$, where $k=2m$, we obtain an integer matrix whose leading principal minors vanish exactly when $\det M_k(X_r)=0$.
These minors were evaluated for all even $k\le 2\ell_0-2$; the integers involved have tens of thousands of digits, and the whole computation for $X_3$ takes about two minutes.

\emph{Step 4: the outcome.}
\begin{itemize}
\item[(a)] $M_k(X_3)=O$ for $k\in S_D$.
\item[(b)] For even $k\le2\ell_0-2$ with $k\notin S_D$, we have $\det M_k(X_r)=0$ exactly for
\[
\begin{array}{c|l}
 X_1 & \text{none}\\
 X_2 & 8,12,16,18,20,26,28,30,34,38,46\\
 X_3 & 8,12,16,18,20,24,26,28,30,32,34,36,38,40,42,44,46,50,52,54,58,62,70
\end{array}
\]
\item[(c)] For each singular $k$ in~(b), the kernel $\ker M_k(X_r)$ contains a vector with all entries positive exactly for
\[
\begin{array}{c|l}
 X_2 & 8,18,26,30,38,46\\
 X_3 & 8,12,18,26,30,34,38,42,46,54,62,70
\end{array}
\]
\end{itemize}
 
By Corollary~\ref{cor:kernel_criterion}, (a)--(c) give the assertions on $\Tev(X_r)$.
For instance, $\ker M_8(X_2)$ is spanned by $(\frac{13}{40},\frac{27}{40})$, and $\ker M_8(X_3)\cap\ker M_{12}(X_3)$ is spanned by $(\frac{97}{960},\frac{243}{560},\frac{625}{1344})$; these recover the weights of the cubature formulas in~\cite[Theorem~$D_4$]{HPV07} (see Theorem~\ref{thm:HPV-D4}~(ii) and~(iii)).

It remains to determine $T(X_r)$.
The equal weight $w\equiv1/|X_r|$ corresponds to the vector $\boldsymbol{c}_r:=(|\mathcal{O}_1|,\ldots,|\mathcal{O}_r|)/|X_r|$, which is proportional to $(1)$, $(1,4)$ and $(1,4,6)$ for $r=1,2,3$, respectively.
By Remark~\ref{rem:h-orbit}, an even integer $k$ belongs to $T(X_r)$ if and only if $\boldsymbol{c}_r\in\ker M_k(X_r)$.
This holds for every $k\in S_D$ by~(a).
Conversely, since $T(X_r)\subset T_\ast(X_r)$, any other even $k\in T(X_r)$ would be one of the integers listed in~(c), and for each of them a direct computation shows $\boldsymbol{c}_r\,M_k(X_r)\ne\boldsymbol{0}$.
Hence $T(X_r)\cap2\N=S_D$.
Finally, since $X_r$ is antipodal, every odd integer belongs to $T(X_r)$, and hence to $T_\ast(X_r)$.
\end{proof}
 
\begin{remark}\label{rem:D4-invariants}
Part~(a) and part of~(b) can be explained without computation.
For $P\in\mathrm{Harm}_k$, the average $|\mathcal{O}_i|^{-1}\sum_{\x\in\mathcal{O}_i}P(\x)$ depends only on the $H$-invariant part of $P$, so that the rank of $M_k(X_3)$ is at most $d_k:=\dim\mathrm{Harm}_k^H$.
By~\cite[Section~6.4]{HPV07},
\[
 \sum_{k\ge0}d_k\,t^k=\frac{1+t^{18}}{(1-t^8)(1-t^{12})(1-t^{24})},
\]
so that, among the even $k>0$, one has $d_k=0$ exactly for $k\in S_D$ and $d_k=1$ exactly for $k\in\{8,12,16,18,20,26,28,30,34,38,46\}$.
For these eleven values of $k$, the matrices $M_k(X_2)$ and $M_k(X_3)$ have rank at most one and are therefore singular; in fact, they are exactly the singular indices of $X_2$ outside $S_D$.
The remaining singular indices $24,32,36,40,42,44,50,52,54,58,62,70$ of $X_3$ are not explained in this way, since $d_k\ge3$ for these $k$.
\end{remark}

\subsection{Shells of $E_8$}\label{app:E8}
Let
\[
 E_8:=\Bigl\{\x\in\Z^8\cup\bigl(\Z+\tfrac12\bigr)^8\ \Bigm|\ x_1+\cdots+x_8\in2\Z\Bigr\}
\]
be the root lattice of type $E_8$, which is the unique even unimodular lattice of rank $8$ up to isometry, and put $\Lambda:=E_8$.
The number of vectors of $\Lambda$ of squared norm $2m$ is $240\,\sigma_3(m)$, where $\sigma_3(m):=\sum_{e\mid m}e^3$.
In particular, its shells of squared norms $2$, $4$, $6$ and $8$ consist of $240$, $2160$, $6720$ and $17520$ vectors, respectively.
Let $H:=\operatorname{Aut}(E_8)$, which is the Weyl group $W(E_8)$.

The group $H$ preserves every normalized shell of $\Lambda$, and $\Lambda[2]$, $\Lambda[4]$ and $\Lambda[6]$ are single $H$-orbits.
The normalized shell $\Lambda[8]$ is not: since $2\x$ has squared norm $8$ for every $\x\in\Lambda$ of squared norm $2$, we have $\Lambda[2]\subset\Lambda[8]$, and $\Lambda[8]\setminus\Lambda[2]$ is a single $H$-orbit.
Accordingly, we put
\[
 \mathcal{O}_1:=\Lambda[2],\qquad
 \mathcal{O}_2:=\Lambda[4],\qquad
 \mathcal{O}_3:=\Lambda[6],\qquad
 \mathcal{O}_4:=\Lambda[8]\setminus\Lambda[2],
\]
so that these four orbits are pairwise disjoint and
\[
 |\mathcal{O}_1|=240,\qquad
 |\mathcal{O}_2|=2160,\qquad
 |\mathcal{O}_3|=6720,\qquad
 |\mathcal{O}_4|=17280.
\]

Every normalized shell of $\Lambda$ is a spherical $7$-design, and $\Lambda[2m]$ is a spherical $8$-design if and only if Ramanujan's $\tau(m)$ vanishes (see~\cite{Ven01} and~\cite[Section~1]{HPV07}).
Combining several shells, de la Harpe, Pache and Venkov obtained the following.

\begin{theorem}[{\cite[Theorem~$E_8$]{HPV07}}]\label{thm:HPV-E8}
\begin{itemize}
\item[$(\mathrm{i})$] $Y_2:=\Lambda[2,4]$ is a weighted $11$-design with $|Y_2|=2400$.
\item[$(\mathrm{ii})$] $Y_3:=\Lambda[2,6,8]$ is a weighted $13$-design with $|Y_3|=24240$.
\item[$(\mathrm{iii})$] $Y_4:=\Lambda[2,4,6,8]$ is a weighted $15$-design with $|Y_4|=26400$.
\end{itemize}
\end{theorem}

Assertion~(i) appears already in~\cite{GS81}.
The orbit decompositions of these three sets are
\[
 Y_2=\mathcal{O}_1\sqcup\mathcal{O}_2,\qquad
 Y_3=\mathcal{O}_1\sqcup\mathcal{O}_3\sqcup\mathcal{O}_4,\qquad
 Y_4=\mathcal{O}_1\sqcup\mathcal{O}_2\sqcup\mathcal{O}_3\sqcup\mathcal{O}_4,
\]
and we always list the orbits in this order.
As observed in~\cite[Section~6.2]{HPV07}, the analogous construction on $\Lambda[2,4,6]$ or $\Lambda[2,4,8]$ leads to a negative weight, which is why the set in~(ii) consists of the shells of squared norms $2$, $6$ and $8$.

Applying the procedure of the proof of Theorem~\ref{thm:D4-shells} with $d=8$ and $H=W(E_8)$, we obtain the following complete description of $T(Y_r)$ and $T_\ast(Y_r)$ for $r=2,3,4$.
The computation is parallel to that in the proof of Theorem~\ref{thm:D4-shells}; we omit the details and record below the data needed to reproduce it.

\begin{theorem}\label{thm:E8-shells}
Let $Y_2$, $Y_3$ and $Y_4$ be as in Theorem~\ref{thm:HPV-E8}, and put
\[
 S_E:=\{2,4,6,10\}.
\]
Then the Laplace cutoffs of $Y_2$, $Y_3$ and $Y_4$ are $30$, $122$ and $134$, respectively,
\[
 T(Y_2)=T(Y_3)=T(Y_4)=S_E\cup\mathbb{O},
\]
and
\begin{align*}
 \Tev(Y_2)&=S_E\cup\{8,14,18,22\},\\
 \Tev(Y_3)&=\Tev(Y_2)\cup\{12,16,26\},\\
 \Tev(Y_4)&=\Tev(Y_3)\cup\{30,34,38,42,46\}.
\end{align*}
\end{theorem}

With $m_{ij}(b)$ defined as in Step~1 of the proof of Theorem~\ref{thm:D4-shells}, the inner product distributions between the orbits $\mathcal{O}_1,\ldots,\mathcal{O}_4$ are listed in Table~\ref{tab:E8-dist}; the remaining ones follow from $|\mathcal{O}_i|\,m_{ij}(b)=|\mathcal{O}_j|\,m_{ji}(b)$.
In particular, $a_{Y_2}^2=\frac{9}{16}$, $a_{Y_3}^2=\frac{49}{64}$ and $a_{Y_4}^2=\frac{25}{32}$.
For even $k$, the $(i,j)$-entry of $M_k(Y_4)$ is $\frac{h_k}{|\mathcal{O}_j|}\sum_{b}m_{ij}(b)\,G_{k,8}(\sqrt b)$, and $M_k(Y_2)$ and $M_k(Y_3)$ are the principal submatrices of $M_k(Y_4)$ indexed by $\{1,2\}$ and $\{1,3,4\}$, respectively.

\begin{table}[ht]
\centering
\small
$\begin{array}{c|l}
 (i,j) & b\ \bigl(m_{ij}(b)\bigr)\\ \hline
 (1,1) & 0\,(126),\ \tfrac{1}{4}\,(112),\ 1\,(2)\\[2pt]
 (1,2) & 0\,(756),\ \tfrac{1}{8}\,(1152),\ \tfrac{1}{2}\,(252)\\[2pt]
 (1,3) & 0\,(2072),\ \tfrac{1}{12}\,(3024),\ \tfrac{1}{3}\,(1512),\ \tfrac{3}{4}\,(112)\\[2pt]
 (1,4) & 0\,(4032),\ \tfrac{1}{16}\,(8064),\ \tfrac{1}{4}\,(4032),\ \tfrac{9}{16}\,(1152)\\[2pt]
 (2,2) & 0\,(574),\ \tfrac{1}{16}\,(896),\ \tfrac{1}{4}\,(560),\ \tfrac{9}{16}\,(128),\ 1\,(2)\\[2pt]
 (2,3) & 0\,(1288),\ \tfrac{1}{24}\,(2688),\ \tfrac{1}{6}\,(1680),\ \tfrac{3}{8}\,(896),\ \tfrac{2}{3}\,(168)\\[2pt]
 (2,4) & 0\,(3360),\ \tfrac{1}{32}\,(5376),\ \tfrac{1}{8}\,(4608),\ \tfrac{9}{32}\,(2688),\ \tfrac{1}{2}\,(1120),\ \tfrac{25}{32}\,(128)\\[2pt]
 (3,3) & 0\,(1260),\ \tfrac{1}{36}\,(1944),\ \tfrac{1}{9}\,(1782),\ \tfrac{1}{4}\,(1084),\ \tfrac{4}{9}\,(540),\ \tfrac{25}{36}\,(108),\ 1\,(2)\\[2pt]
 (3,4) & 0\,(2304),\ \tfrac{1}{48}\,(5184),\ \tfrac{1}{12}\,(3888),\ \tfrac{3}{16}\,(3168),\ \tfrac{1}{3}\,(1728),\ \tfrac{25}{48}\,(864),\ \tfrac{3}{4}\,(144)\\[2pt]
 (4,4) & 0\,(2310),\ \tfrac{1}{64}\,(4048),\ \tfrac{1}{16}\,(3920),\ \tfrac{9}{64}\,(2800),\ \tfrac{1}{4}\,(2296),\ \tfrac{25}{64}\,(1232),\\
       & \tfrac{9}{16}\,(560),\ \tfrac{49}{64}\,(112),\ 1\,(2)
\end{array}$\vspace{5mm}
\caption{The squared inner products $b$ and the numbers $m_{ij}(b)$ for $\Lambda=E_8$.}
\label{tab:E8-dist}
\end{table}

\section*{Acknowledgments}
This work was partially supported by
JSPS KAKENHI Grant Numbers 22K03402, 24K06688, 24K06871, 26K00617, 26K00620, 26K22252 and 26K22262, and the Research Institute for Mathematical Sciences, an International Joint Usage/Research Center located in Kyoto University.

\section*{Statements and Declarations}
We wish to confirm that there are no known conflicts of interest associated with this publication and there has been no significant financial support for this work that could have influenced its outcome.

The data that support the findings of this study are available from the corresponding author, K.T., upon reasonable request.

\section*{Use of generative AI}
During the preparation of this work, the authors used Claude (Anthropic)
to improve the language and presentation of the manuscript, to check
notation, references and typographical errors, and to independently
verify the computations in Section~7 and Appendix~A. The authors
reviewed and edited the content as needed and take full responsibility
for the content of the article.



\begin{thebibliography}{99}

\bibitem{And00}
N.~N.~Andreev,
\textit{A minimal design of order $11$ on the $3$-sphere},
Math. Notes \textbf{67} (2000), no.~4, 417--424;
translated from Mat.\ Zametki \textbf{67} (2000), no.~4, 489--497.

\bibitem{Askey}
R.~Askey,
\textit{Orthogonal Polynomials and Special Functions},
CBMS-NSF Regional Conf.\ Ser.\ in Appl.\ Math.\ \textbf{21},
SIAM, Philadelphia, PA, 1975.

\bibitem{Baj07}
B.~Bajnok,
\textit{Orbits of the hyperoctahedral group as Euclidean designs},
J.\ Algebraic Combin.\ \textbf{25} (2007), no.~4, 375--397.

\bibitem{BB09}
E.~Bannai and E.~Bannai,
\textit{A survey on spherical designs and algebraic combinatorics on spheres},
European J.\ Combin.\ \textbf{30} (2009), no.~6, 1392--1425.

\bibitem{BBHS10}
E.~Bannai, E.~Bannai, M.~Hirao, and M.~Sawa,
\textit{Cubature formulas in numerical analysis and {E}uclidean tight designs},
European J.\ Combin.\ \textbf{31} (2010), no.~2, 423--441.

\bibitem{BOT15}
E.~Bannai, T.~Okuda, and M.~Tagami,
\textit{Spherical designs of harmonic index $t$},
J.\ Approx.\ Theory \textbf{195} (2015), 1--18.

\bibitem{BDK01}
P.~Boyvalenkov, D.~Danev, and P.~Kazakov,
\textit{Indexes of spherical codes},
in \textit{Codes and Association Schemes},
DIMACS Ser.\ Discrete Math.\ Theoret.\ Comput.\ Sci.\ \textbf{56},
Amer.\ Math.\ Soc., Providence, RI, 2001, pp.~47--57.

\bibitem{Del73}
P.~Delsarte,
\textit{An Algebraic Approach to the Association Schemes of Coding Theory},
Philips Res.\ Repts.\ Suppl.\ \textbf{10} (1973).

\bibitem{DGS77}
P.~Delsarte, J.~M.~Goethals, and J.~J.~Seidel,
\textit{Spherical codes and designs},
Geom.\ Dedicata \textbf{6} (1977), 363--388.

\bibitem{DL98}
P.~Delsarte and V.~I.~Levenshtein,
\textit{Association schemes and coding theory},
IEEE Trans.\ Inform.\ Theory \textbf{44} (1998), no.~6, 2477--2504.

\bibitem{DS89}
P.~Delsarte and J.~J.~Seidel,
\textit{Fisher type inequalities for Euclidean $t$-designs},
Linear Algebra Appl.\ \textbf{114/115} (1989), 213--230.

\bibitem{GS81}
J.~M.~Goethals and J.~J.~Seidel,
\textit{Cubature formulae, polytopes, and spherical designs},
in \textit{The Geometric Vein: The Coxeter Festschrift},
C.~Davis, B.~Gr\"{u}nbaum, and F.~A.~Sherk (eds.),
Springer, New York--Berlin, 1981, pp.~203--218.

\bibitem{HS96}
R.~H.~Hardin and N.~J.~A.~Sloane,
\textit{McLaren's improved snub cube and other new spherical designs
in three dimensions},
Discrete Comput.\ Geom.\ \textbf{15} (1996), no.~4, 429--441.


\bibitem{HPV07}
P.~de~la~Harpe, C.~Pache, and B.~Venkov,
\textit{Construction of spherical cubature formulas using lattices},
St.\ Petersburg Math.\ J.\ \textbf{18} (2007), no.~1, 119--139.

\bibitem{HNT23}
M.~Hirao, H.~Nozaki, and K.~Tasaka,
\textit{Spherical designs and modular forms of the $D_4$ lattice},
Res.\ Number Theory \textbf{9} (2023), no.~4, Paper No.~77, 18~pp.

\bibitem{HNT}
M.~Hirao, H.~Nozaki, and K.~Tasaka,
\textit{Spherical designs for finite quaternionic unit groups
and their applications to modular forms},
Algebr.\ Comb.\ \textbf{9} (2026), no.~4, 1157--1178.

\bibitem{Kra17}
I.~Krasikov,
\textit{On approximation of ultraspherical polynomials in the oscillatory region},
J.\ Approx.\ Theory \textbf{222} (2017), 143--156.

\bibitem{Mie13}
T.~Miezaki,
\textit{On a generalization of spherical designs},
Discrete Math.\ \textbf{313} (2013), no.~4, 375--380.

\bibitem{MN26}
R.~Misawa and Y.~Nishimura,
\textit{Spherical designs with infinite harmonic strength},
arXiv:2607.01761, 2026.

\bibitem{Mol76}
H.~M.~M\"{o}ller,
\textit{Kubaturformeln mit minimaler Knotenzahl},
Numer.\ Math.\ \textbf{25} (1976), 185--200.

\bibitem{Mol79}
H.~M.~M\"{o}ller,
\textit{Lower bounds for the number of nodes in cubature formulae},
in \textit{Numerische Integration},
G.~H\"{a}mmerlin (ed.),
Internat.\ Ser.\ Numer.\ Math.\ \textbf{45},
Birkh\"{a}user, Basel, 1979, pp.~221--230.

\bibitem{NS88}
A.~Neumaier and J.~J.~Seidel,
\textit{Discrete measures for spherical designs, eutactic stars and lattices},
Indag.\ Math.\ \textbf{50} (1988), 321--334.

\bibitem{NS12}
H.~Nozaki and M.~Sawa,
\textit{Note on cubature formulae and designs obtained from group orbits},
Canad. J. Math. \textbf{64} (2012), no.~6, 1359--1377.


\bibitem{OY16}
T.~Okuda and W.-H.~Yu,
\textit{A new relative bound for equiangular lines and nonexistence of
tight spherical designs of harmonic index $4$},
European J.\ Combin.\ \textbf{53} (2016), 96--103.


\bibitem{Pan22}
B.~V.~Pandey,
\textit{Modular forms and ellipsoidal $T$-designs},
Ramanujan J.\ \textbf{58} (2022), 1245--1257.

\bibitem{Sob62}
S.~L.~Sobolev,
\textit{Cubature formulas on the sphere invariant under finite groups of rotations},
Dokl.\ Akad.\ Nauk SSSR \textbf{146} (1962), 310--313;
English translation: Soviet Math.\ Dokl.\ \textbf{3} (1962), 1307--1310.

\bibitem{Sze75}
G.~Szeg\H{o},
\textit{Orthogonal Polynomials},
4th ed., Amer.\ Math.\ Soc.\ Colloq.\ Publ.\ \textbf{23},
American Mathematical Society, Providence, RI, 1975.

\bibitem{Ven01}
B.~Venkov,
\textit{R\'{e}seaux et designs sph\'{e}riques},
in \textit{R\'{e}seaux euclidiens, designs sph\'{e}riques et formes modulaires},
J.~Martinet (ed.),
Monogr.\ Enseign.\ Math.\ \textbf{37},
Enseignement Math., Geneva, 2001, pp.~10--86.

\bibitem{ZBBKY17}
Y.~Zhu, E.~Bannai, E.~Bannai, K.-T.~Kim, and W.-H.~Yu,
\textit{On spherical designs of some harmonic indices},
Electron.\ J.\ Combin.\ \textbf{24} (2017), no.~2, Paper No.~2.14, 28~pp.

\end{thebibliography}
\end{document}